\magnification=1200
\loadmsam
\loadmsbm
\loadeufm
\loadeusm
\UseAMSsymbols     
\input amssym.def

\font\BIGtitle=cmr10 scaled\magstep3
\font\bigtitle=cmr10 scaled\magstep1
\font\boldsectionfont=cmb10 scaled\magstep1
\font\section=cmsy10 scaled\magstep1

\def\scr#1{{\fam\eusmfam\relax#1}}
\def\scrA{{\scr A}}
\def\scrB{{\scr B}}
\def\scrC{{\scr C}}
\def\scrD{{\scr D}}
\def\scrE{{\scr E}}
\def\scrF{{\scr F}}
\def\scrG{{\scr G}}

\def\scrH{{\scr H}}
\def\scrI{{\scr I}}
\def\scrL{{\scr L}}
\def\scrK{{\scr K}}
\def\scrJ{{\scr J}}
\def\scrM{{\scr M}}
\def\scrN{{\scr N}}
\def\scrO{{\scr O}}
\def\scrP{{\scr P}}

\def\scrS{{\scr S}}
\def\scrU{{\scr U}}
\def\scrR{{\scr R}}

\def\scrV{{\scr V}}
\def\scrX{{\scr X}}
\def\scrY{{\scr Y}}
\def\scrW{{\scr W}}

\def\gr#1{{\fam\eufmfam\relax#1}}

\def\grC{{\gr C}}	
\def\grD{{\gr D}}

	\def\gri{{\gr i}}
	
\def\grK{{\gr K}}	\def\grk{{\gr k}}

	\def\grr{{\gr r}}
\def\grS{{\gr S}}	\def\grs{{\gr s}}

	\def\grw{{\gr w}}

\def\db#1{{\fam\msbfam\relax#1}}

\def\dbC{{\db C}} \def\dbD{{\db D}}
 \def\dbF{{\db F}}
\def\dbG{{\db G}}

 \def\dbN{{\db N}}
 
\def\dbQ{{\db Q}} \def\dbR{{\db R}}
\def\dbS{{\db S}} \def\dbT{{\db T}}
 
\def\dbW{{\db W}} 
 \def\dbZ{{\db Z}}

\def\eps{{\varepsilon}}

\def\im{\text{Im}}

\def\Ker{\text{Ker}}

\def\Sh{\hbox{\rm Sh}}

\def\red{\text{red}}

\def\ad{\text{ad}}

\def\Hom{\text{Hom}}
\def\End{\text{End}}
\def\Spec{\text{Spec}}

\def\red{\text{red}}

\def\Lie{\text{Lie}}

\def\leaderfill{\leaders\hbox to 1em
     {\hss.\hss}\hfill}
\def\nspace{\lineskip=1pt\baselineskip=12pt\lineskiplimit=0pt}

\def\finishproclaim{\par\rm
     \ifdim\lastskip<\medskipamount\removelastskip
     \penalty55\medskip\fi}
\def\endproof{$\hfill \square$}
\def\proof{\par\noindent {\it Proof:}\enspace}
\def\references#1{\par
  \centerline{\boldsectionfont References}\medskip
     \parindent=#1pt\nspace}
\def\Ref[#1]{\par\hang\indent\llap{\hbox to\parindent
     {[#1]\hfil\enspace}}\ignorespaces}
\def\Item#1{\par\smallskip\hang\indent\llap{\hbox to\parindent
     {#1\hfill$\,\,$}}\ignorespaces}
\def\ItemItem#1{\par\indent\hangindent2\parindent
     \hbox to \parindent{#1\hfill\enspace}\ignorespaces}

\def\arrowsim{\,\smash{\mathop{\to}\limits^{\lower1.5pt
  \hbox{$\scriptstyle\sim$}}}\,}

\def\doublemaprights#1#2#3#4{\raise3pt\hbox{$\mathop{\,\,\hbox to     
#1pt{\rightarrowfill}\kern-30pt\lower3.95pt\hbox to
     #2pt{\rightarrowfill}\,\,}\limits_{#3}^{#4}$}}

\def\rightcapdownarrow{\raise9pt\hbox{$\ssize\cap$}\kern-7.75pt
     \Big\downarrow}

\def\rcapmapdown#1{\rightcapdownarrow\kern-1.0pt\vcenter{
     \hbox{$\scriptstyle#1$}}}

\def\rmapdown#1{\Big\downarrow\kern-1.0pt\vcenter{
     \hbox{$\scriptstyle#1$}}}
\def\rightsubsetarrow#1{{\ssize\subset}\kern-4.5pt\lower2.85pt
     \hbox to #1pt{\rightarrowfill}}
\def\longtwoheadedrightarrow#1{\raise2.2pt\hbox to #1pt{\hrulefill}
     \!\!\!\twoheadrightarrow}

\NoBlackBoxes
\parindent=25pt
\document
\footline={\hfil}
\footline={\hss\tenrm \folio\hss}
\pageno=1
\bigskip

\bigskip
\noindent 
\centerline{\BIGtitle Level m stratifications}
\bigskip
\centerline{\BIGtitle  of versal deformations of p-divisible groups}
\bigskip
\centerline{\bigtitle Adrian Vasiu, Binghamton University}
\medskip
\centerline{accepted (in final form) for publication in J. Alg. Geom.}
\medskip
\centerline{\it To Carlo Traverso, for his 62nd birthday} 

\bigskip\noindent
{\bf ABSTRACT.} Let $k$ be an algebraically closed field of characteristic $p>0$. Let $c,d,m$ be positive integers. Let $D$ be a $p$-divisible group of codimension $c$ and dimension $d$ over $k$. Let $\scrD$ be a versal deformation of $D$ over a smooth $k$-scheme $\scrA$ which is equidimensional of dimension $cd$. We show that there exists a reduced, locally closed subscheme $\grs_D(m)$ of $\scrA$ that has the following property: a point $y\in\scrA(k)$ belongs to $\grs_D(m)(k)$ if and only if $y^*(\scrD)[p^m]$ is isomorphic to $D[p^m]$. We prove that $\grs_D(m)$ is {\it regular and equidimensional} of {\it dimension} $cd-\dim(\pmb{\text{Aut}}(D[p^m]))$. We give a proof of {\it Traverso's formula} which for $m>>0$ computes the codimension of $\grs_D(m)$ in $\scrA$ (i.e., $\dim(\pmb{\text{Aut}}(D[p^m]))$) in terms of the Newton polygon of $D$. We also provide a criterion of when $\grs_D(m)$ satisfies the {\it purity property} (i.e., it is an affine $\scrA$-scheme). Similar results are proved for {\it quasi Shimura $p$-varieties of Hodge type} that generalize the special fibres of good integral models of Shimura varieties of Hodge type in unramified mixed characteristic $(0,p)$.



\bigskip\bigskip
\noindent
{\boldsectionfont 1. Introduction}
\bigskip

Let $p\in\dbN$ be a prime. Let $k$ be an algebraically closed field of characteristic $p$. Let $c$ and $d$ be positive integers and let $r:=c+d$. Let $D$ be a $p$-divisible group of codimension $c$ and dimension $d$ over $k$. The height of $D$ is $r$. Let $n_D$ be the smallest positive integer for which the following statement holds: if $D_1$ is a $p$-divisible group over $k$ such that $D_1[p^{n_D}]$ is isomorphic to $D[p^{n_D}]$, then $D_1$ is isomorphic to $D$. For the existence of $n_D$ we refer to [Ma, Ch. III, Sect. 3], [Tr1, Thm. 3], [Tr2, Thm. 1], [Va2, Cor. 1.3], or [Oo2, Cor. 1.7]. 

Let $W(k)$ be the ring of Witt vectors
 with coefficients in $k$. Let $B(k)$ be the field of fractions of $W(k)$. Let $\sigma:=\sigma_k$ be the Frobenius automorphism of $W(k)$ and $B(k)$ induced from $k$. Let $(M,\phi)$ be the (contravariant) Dieudonn\'e module of $D$. We recall that $M$ is a free $W(k)$-module of rank $r$ and $\phi:M\to M$ is a $\sigma$-linear endomorphism such that we have $pM\subseteq \phi(M)$. Dieudonn\'e's classification of $F$-isocrystals over $k$ (see [Di, Thms. 1 and 2], [Ma, Ch. 2, Sect. 4], [Dem, Ch. IV], etc.) implies that we have a direct sum decomposition $(M[{1\over p}],\phi)=\oplus_{s=1}^v (W_s,\phi)$
into simple $F$-isocrystals over $k$ (here $v$ is a positive integer). More precisely, for $s\in\{1,\ldots,v\}$ there exist $c_s,d_s\in\dbN\cup\{0\}$ such that $r_s:=c_s+d_s>0$, $g.c.d.\{c_s,d_s\}=1$, $\dim_{B(k)}(W_s)=r_s$, and moreover there exists a $B(k)$-basis for $W_s$ formed by elements fixed by $p^{-d_s}\phi^{r_s}$; the unique Newton polygon slope of $(W_s,\phi)$ is $\alpha_s:={{d_s}\over {r_s}}\in \dbQ\cap [0,1]$. 

Let $\scrA$ be a smooth $k$-scheme which is equidimensional of dimension $cd$ and for which the following two properties hold:

\medskip
{\bf (i)} there exists a $p$-divisible group $\scrD$ of codimension $c$ and dimension $d$ over $\scrA$  which is a versal deformation at each $k$-valued point of $\scrA$;

\smallskip
{\bf (ii)} there exists a point $y_D\in\scrA(k)$ such that $y_D^*(\scrD)$ is isomorphic to $D$.

\medskip
In practice, $\scrA$ is the pull back to $\Spec(k)$ of a special fibre of a good integral model in unramified mixed characteristic $(0,p)$ of a {\it unitary Shimura variety} $\Sh(\scrG,\scrX)$ which is constructed as in [Zi, Subsect. 3.5] or [Ko, Sect. 5], where the real adjoint group $\scrG_{\dbR}^{\ad}$ is isomorphic to $\pmb{PGU}(c,d)\times_{\dbR} \pmb{PGU}(c+d)^u$ for some $u\in\dbN\cup\{0\}$ and where $\scrG_{\dbQ_p}$ is a split group over $\dbQ_p$ (see [Va5, Example 5.4.3 (b)] and [Va6] for more details on the fact that such a triple $(\scrA,\scrD,y_D)$ always exists). 

Let $m$ be a positive integer. The first goal of the paper is to study the following set of $k$-valued points 
$$\grs_D(m)(k):=\{y\in\scrA(k)|y^*(\scrD)[p^m]\;\text{is}\;\text{isomorphic}\;\text{to}\;D[p^m]\}.$$
The second goal of the paper is to study analogues of the set $\grs_D(m)(k)$ that pertain to good integral models in unramified mixed characteristic $(0,p)$ of {\it Shimura varieties of Hodge type} and that define {\it level $m$ stratifications} of special fibres of such integral models. Here and in the whole paper the stratifications of reduced schemes over a field are as defined in [Va2, Subsubsect. 2.1.1]. We recall that Shimura varieties of Hodge type are {\it moduli schemes} of polarized abelian schemes endowed with families of Hodge cycles and with symplectic similitude structures. These {\it new} level $m$ stratifications:

\medskip
{\bf (a)} (for $m=1$) generalize the {\it Ekedahl--Oort stratifications} studied in [Oo1], [Mo], [We], etc., and the {\it mod $p$ stratifications} studied in [Va3, Sect. 12];

\smallskip
{\bf (b)} (for $m>>0$)  generalize the {\it ultimate stratifications} studied in [Va2, Thm. 5.3.1 and Subsubsect. 5.3.2] and the {\it foliations} studied in [Oo2];

\smallskip
{\bf (c)} represent a fundamental tool in studying good integral models of Shimura varieties of Hodge type (like their special fibres on which groups of Hecke orbits act, their cohomology groups, their local and global geometries, their crystalline properties, etc.).

\medskip
In order to state our Basic Theorem, we will need the following two definitions that recall [Va4, Def. 1.6.1] and the very essence of the purity property introduced in [Va2, Subsubsect. 2.1.1].

\medskip\smallskip\noindent
{\bf 1.1. Definition.} {\bf (a)} Let $\pmb{\text{Aut}}(D[p^m])$ be the group scheme over $k$ of automorphisms of $D[p^m]$ and let $\gamma_D(m):=\dim(\pmb{\text{Aut}}(D[p^m]))$. We call $(\gamma_D(m))_{m\ge 1}$ the {\it centralizing sequence} of $D$. We also refer to $s_D:=\gamma_D(n_D)$ as the {\it specializing height} of $D$. 

\smallskip
{\bf (b)} A reduced, locally closed subscheme of a reduced $k$-scheme $\scrM$ is said to satisfy the {\it purity property}, if it is an affine $\scrM$-scheme.   

\medskip\smallskip\noindent
{\bf 1.2. Basic Theorem.} {\it With the above notations, the following eight properties hold:

\medskip
{\bf (a)} there exists a unique reduced, locally closed subscheme $\grs_D(m)$ of $\scrA$ such that our notations match i.e., the set $\grs_D(m)(k)$ we introduced above is the set of $k$-valued points of $\grs_D(m)$;

\smallskip
{\bf (b)} the scheme $\grs_D(m)$ is regular and equidimensional;

\smallskip
{\bf (c)} we have $\dim(\grs_D(m))=cd-\gamma_D(m)$;

\smallskip
{\bf (d)} for $m\ge n_D$, we have $\grs_D(m)=\grs_D(n_D)$ and $\gamma_D(m)=s_D$;

\smallskip
{\bf (e)} the specializing height $s_D$ of $D$ is an isogeny invariant;

\smallskip
{\bf (f)} we have $\dim(\grs_D(n_D))={1\over 2}\sum_{s=1}^v\sum_{t=1}^v r_sr_t|\alpha_s-\alpha_t|={1\over 2}\sum_{s=1}^v\sum_{t=1}^v |c_sd_t-c_td_s|$;

\smallskip
{\bf (g)} the $\scrA$-scheme $\grs_D(m)$ is quasi-affine;

\smallskip
{\bf (h)} if $m\ge 2$ and if the image of the homomorphism $\pmb{\text{Aut}}(D[p^m])(k)\to\pmb{\text{Aut}}(D[p])(k)$ is finite, then the reduced, locally closed subscheme $\grs_D(m)$ of $\scrA$ satisfies the {\it purity property}.}

\medskip
For the sake of completeness we also state here the following practical property that was conjectured by Traverso (cf. [Tr3, Sect. 40, Conj. 5]) and that is proved in [NV, Thm. 1.2]:

\medskip
{\bf (i)} the isogeny class (i.e., the Newton polygon) of $D$ depends only on $D[p^{\lceil{{cd}\over {c+d}}\rceil}]$.

\medskip
For $m\ge 1$ we have a natural monomorphism $\grs_D(m+1)\hookrightarrow\grs_D(m)$. Thus the sequence $\dim(\grs_D(m))_{m\ge 1}$ is decreasing. From this and Theorem 1.2 (c) and (d) we get:

\medskip\noindent
{\bf 1.2.1. Corollary.} {\it For each $p$-divisible group $D$ of codimension $c$ and dimension $d$ over $k$ and for each positive integer $m$, we have $\gamma_D(m)\in\{0,1,\ldots,s_D\}\subseteq\{0,1,\ldots,cd\}$.}

\medskip\smallskip\noindent
{\bf 1.3. On literature.} The fact that $\grs_D(m)(k)$ is a constructible subset of $\scrA(k)$ is a standard piece of algebraic geometry. The classification of commutative, finite group schemes over $k$ annihilated by $p$ accomplished by Kraft in [Kr], implies that $\scrA$ is a finite disjoint union of reduced, locally closed subschemes of it of the form $\grs_D(1)$. Theorem 1.2 (b) is a direct consequence of Grothendieck's results on local deformations of truncated Barsotti--Tate groups presented in [Il]. The facts that $\grs_D(1)$ is an equidimensional $k$-scheme and a quasi-affine $\scrA$-scheme are only  a variant of [Oo1, Thm. (1.2)]. The fact that $\grs_D(m)$ is a regular and equidimensional $k$-scheme is first proved in [Va2, Basic Thm. 5.3.1 (b) and Rm. 5.3.4 (b)]; a variant of this for $m=n_D$ also shows up in [Oo2, Thm. 3.13]. The formula $\dim(\grs_D(1))=cd-\gamma_D(1)$ is proved first for $p>2$ in [We, Thm. of Introd.] (see also [MW, Subsects. 7.10 to 7.14]) and for all $p$ in [Va3, Basic Thms. A and D]. For $m\ge 2$, Theorem 1.2 (c) was not previously available in the literature; however, it can be proved as well using the ideas of the proof of [We, Thm. of Introd.]. Theorem 1.2 (d) is a consequence of the definition of $n_D$. A variant of Theorem 1.2 (e) was first obtained by Traverso, cf. [Tr2, Thm. 2]. Theorem 1.2 (f) was first obtained by Traverso (cf. [Tr2, Sect. 1, p. 48]) but it never got published; it also shows up in an informal manuscript of Oort. For $m=n_D$, Theorem 1.2 (h) is implied by [Va2, Thm. 5.3.1 (c)]. For $2\le m\le n_D-1$, Theorem 1.2 (h) is new. It seems to us that for $m\ge 2$, the Corollary 1.2.1 is new. 

\medskip\smallskip\noindent
{\bf 1.4. On contents.} In Section 2 we follow [Va4, Sect. 5] in order to introduce {\it orbit spaces} of {\it truncated Barsotti--Tate groups of level $m$} over $k$ that have codimension $c$ and dimension $d$. Such orbit spaces were first considered in [Tr2, Thm. 2] and [Tr3, Sects. 26 to 40]. A significant change in their presentation was made independently of [Tr2] and [Tr3], in [Va3] for the case $m=1$ and in [Va4, Sect. 5] for all $m$. The change allows a very easy description of the orbit spaces that leads to short, elementary, and foundational (computations and) proofs of all parts of the Basic Theorem; the change is explained in Remark 2.4.1. The proof of the Basic Theorem is presented in Section 3. In Section 4 we show how the Basic Theorem gets easily translated to the case of {\it quasi Shimura $p$-varieties of Hodge type} that generalize the special fibres of good integral models in unramified mixed characteristic of Shimura varieties of Hodge type. In particular, the Basic Theorem gets easily translated to the case of special fibres of Mumford's moduli schemes $\scrA_{d,1,l}$ (see Example 4.5). For $m=1$, the Basic Corollary 4.3 (b) generalizes results of [We] and [Mo] obtained for special fibres of good integral models of Shimura varieties of PEL type.

\bigskip\noindent
{\boldsectionfont 2. Orbit spaces of truncated Barsotti--Tate groups of level m}
\bigskip

In this Section we recall the group action $\dbT_m$ over $k$ we introduced in [Va4, Sect. 5]; its set of orbits parametrizes isomorphism classes of truncated Barsotti--Tate groups of level $m$ over $k$ that have codimension $c$ and dimension $d$. Subsection 2.1 introduces certain group schemes that play a key role in the definition (see Section 2.2) of $\dbT_m$. Subsection 2.3 studies $\dbT_1$. Theorem 2.4 recalls properties of stabilizer subgroup schemes of $\dbT_m$ we obtained in [Va4, Thm. 5.3]. Lemmas 2.5 and 2.6 are the very essence of Theorem 1.2 (h).

The notations $p$, $k$, $c$, $d$, $r$, $D$, $n_D$, $W(k)$, $B(k)$, $\sigma$, $(M,\phi)$, $(M[{1\over p}],\phi)=\oplus_{s=1}^v (W_s,\phi)$, $c_s,d_s$, $r_s$, $\alpha_s$, $m$, $\scrA$, $\scrD$, $y_D\in\scrA(k)$, $\grs_D(m)(k)$, $\pmb{\text{Aut}}(D[p^m])$, $\gamma_D(m)$, and $s_D$ are as in Section 1. Let $\vartheta:=p\phi^{-1}:M\to M$ be the Verschiebung map of $(M,\phi)$.

For a commutative $k$-algebra $R$, let $W_m(R)$ be the ring of Witt vectors of length $m$ with coefficients in $R$, let $W(R)$ be the ring of Witt vectors with coefficients in $R$, and let $\sigma_R$ be the Frobenius endomorphism of either $W_m(R)$ or $W(R)$. Let $\delta_m$ be the natural divided power structure on the kernel of the reduction $W(k)$-epimorphism $W_m(R)\twoheadrightarrow R$. Let $\grS_m(R)$ be the thickening $(\Spec(R)\hookrightarrow \Spec(W_m(R)),\delta_m)$ of the Berthelot crystalline site $CRIS(\Spec(R)/\Spec(W(k)))$ introduced in [Be, Ch. III, Sect. 4]. We refer to [BBM] for the crystalline contravariant Dieudonn\'e functor $\dbD$ defined on the category of $p$-divisible groups over $\Spec(R)$. We denote also by $\phi$ the $\sigma$-linear automorphism of $\End(M)[{1\over p}]$ which takes $e\in \End(M)[{1\over p}]$ to $\phi(e):=\phi\circ e\circ\phi^{-1}$.

\medskip\smallskip\noindent
{\bf 2.1. Group schemes.} Let $M=F^1\oplus F^0$ be a direct sum decomposition such that $\bar F^1:=F^1/pF^1$ is the kernel of the reduction  modulo $p$ of $\phi$. Let $\bar F^0:=F^0/pF^0$. The ranks of $F^1$ and $F^0$ are $d$ and $c$ (respectively). The decomposition $M=F^1\oplus F^0$ gives birth naturally to a direct sum decomposition of $W(k)$-modules
$$\End(M)=\Hom(F^0,F^1)\oplus\End(F^1)\oplus\End(F^0)\oplus\Hom(F^1,F^0).$$ 
Let $\scrW_+$ be the maximal subgroup scheme of $\pmb{GL}_{M}$ that fixes both $F^1$ and $M/F^1$; it is a closed subgroup scheme of $\pmb{GL}_M$ whose Lie algebra is the direct summand $\Hom(F^0,F^1)$ of $\End(M)$ and whose relative dimension is $cd$. Let $\scrW_0:=\pmb{GL}_{F^1}\times_{W(k)} \pmb{GL}_{F^0}$; it is a closed subgroup scheme of $\pmb{GL}_M$ whose Lie algebra is the direct summand $\End(F^1)\oplus\End(F^0)$ of $\End(M)$ and whose relative dimension is $d^2+c^2$. The maximal parabolic subgroup scheme $\scrW_{+0}$ of $\pmb{GL}_{M}$ that normalizes $F^1$ is the semidirect product of $\scrW_+$ and $\scrW_0$. Let $\scrW_-$ be the maximal subgroup scheme of $\pmb{GL}_{M}$ that fixes $F^0$ and $M/F^0$; it is a closed subgroup scheme of $\pmb{GL}_M$ whose Lie algebra is the direct summand $\Hom(F^1,F^0)$ of $\End(M)$ and whose relative dimension is $cd$. The maximal parabolic subgroup scheme $\scrW_{0-}$ of $\pmb{GL}_{M}$ that normalizes $F^0$ is the semidirect product of $\scrW_-$ and $\scrW_0$. If $R$ is a commutative $W(k)$-algebra, then we have
$$\scrW_+(R)=1_{M\otimes_{W(k)} R}+\Hom(F^0,F^1)\otimes_{W(k)} R$$ and 
$$\scrW_-(R)=1_{M\otimes_{W(k)} R}+\Hom(F^1,F^0)\otimes_{W(k)} R.$$ 
These identities imply that the group schemes $\scrW_+$ and $\scrW_-$ are isomorphic to $\dbG_a^{cd}$ over $\Spec(W(k))$; in particular, they are smooth and commutative. Let
$$\scrH:=\scrW_+\times_{W(k)}\scrW_0\times_{W(k)} \scrW_-;$$ 
it is a smooth, affine scheme  over $\Spec(W(k))$ of relative dimension $cd+d^2+c^2+cd=r^2$. We consider the natural product morphism $\scrP_0:\scrH\to\pmb{GL}_{M}$ and the following morphism $\scrP_-:=1_{\scrW_+}\times 1_{\scrW_0}\times p1_{\scrW_-}:\scrH\to\scrH$. Let 
$$\scrP_{0-}:=\scrP_0\circ\scrP_{-}:\scrH\to\pmb{GL}_M;$$ 
it is a morphism of $\Spec(W(k))$-schemes whose generic fibre is an open embedding of $\Spec(B(k))$-schemes. 

For $g\in\pmb{GL}_M(W(k))$ and $h=(h_1,h_2,h_3)\in\scrH(W(k))$, let $g[m]\in\pmb{GL}_M(W_m(k))$ and $h[m]=(h_1[m],h_2[m],h_3[m])\in\scrH(W_m(k))$ be the reductions  modulo $p^m$ of $g$ and $h$ (respectively). Thus $1_{M/p^mM}=1_M[m]$. Let $\phi_m,\vartheta_m:M/p^mM\to M/p^mM$ be the reductions  modulo $p^m$ of $\phi,\vartheta:M\to M$.

\medskip\noindent
{\bf 2.1.1. The study of $\scrH$.} Let $\tilde\scrH$ be the {\it dilatation} of $\pmb{GL}_M$ centered on the smooth subgroup $\scrW_{+0k}$ of $\pmb{GL}_{M/pM}$ (see [BLR, Ch. 3, 3.2] for dilatations). We recall that if $\pmb{GL}_M=\Spec(R_M)$ and if $I_{+0k}$ is the ideal of $R_M$ that defines $\scrW_{+0k}$, then as a scheme $\tilde\scrH$ is the spectrum of the $R_M$-subalgebra $R_{\tilde\scrH}$ of $R_M[{1\over p}]$ generated by all elements ${*\over p}$ with $*\in I_{+0k}$ (see loc. cit.). It is well known that $\tilde\scrH$ is a smooth, affine group scheme over $\Spec(W(k))$ which is uniquely determined by the following two additional properties (they follow directly from the definition of $R_{\scrH}$; see [BLR, Ch. 3, 3.2, Props. 1, 2, and 3]):

\medskip
{\bf (i)} there exists a homomorphism $\tilde\scrP_{0-}:\tilde\scrH\to \pmb{GL}_M$ whose generic fibre is an isomorphism of $\Spec(B(k))$-schemes;

\smallskip
{\bf (ii)} a morphism $f:X\to \pmb{GL}_M$ of flat $\Spec(W(k))$-schemes factors (uniquely) through $\tilde\scrP_{0-}$ if and only if the morphism $f_k:X_k\to \pmb{GL}_{M/pM}$ factors through $\scrW_{+0k}$.

\medskip
Due to the property (ii), we have $\tilde\scrH(W(k))=\{*\in\pmb{GL}_M(W(k))|*[1]\in\scrW_{+0}(k)\}$ and $\tilde\scrH(W(k)[\eps]/(\eps^2))=\{*\in\pmb{GL}_M(W(k)[\eps]/(\eps^2))|*\;\text{modulo}\; p\;\text{belongs}\;\text{to}\;\scrW_{+0}(k[\eps]/(\eps^2))\}$.

The image of the map $\scrP_{0-}(W(k)):\scrH(W(k))\to\pmb{GL}_M(W(k))$ is the subgroup $\tilde\scrH(W(k))$ of $\pmb{GL}_M(W(k))$. Due to this and the property (ii), one easily gets that $\scrP_{0-}$ factors through $\tilde\scrP_{0-}$ i.e., there exists a unique affine morphism
$$\scrP:\scrH\to\tilde\scrH$$
such that we have $\scrP_{0-}=\tilde\scrP_{0-}\circ\scrP$. 
The resulting maps $\scrP(W(k)):\scrH(W(k))\to\tilde\scrH(W(k))$ and $\scrP(W(k)[\eps]/(\eps^2)):\scrH(W(k)[\eps]/(\eps^2))\to\tilde\scrH(W(k)[\eps]/(\eps^2))$ are bijections, to be viewed in what follows as natural identifications. From this and the fact that $\tilde\scrH$ and $\tilde\scrH$ are smooth, affine $\Spec(W(k))$-schemes, we get that the morphism $\scrP_k:\scrH_k\to\tilde\scrH_k$ of $\Spec(k)$-schemes is an isomorphism. This implies that the morphism 
$$\scrP_{W_m(k)}:\scrH_{W_m(k)}\to\tilde\scrH_{W_m(k)}$$ 
of smooth, affine $\Spec(W_m(k))$-schemes is an isomorphism, to be viewed in what follows as a natural identification. Therefore $\scrH_{W_m(k)}$ has a natural structure of a smooth, affine group scheme over $\Spec(W_m(k))$. Thus the $p$-adic completion of $\scrH$ has a natural structure of a formal group scheme over $\text{Spf}(W(k))$ isomorphic to the $p$-adic completion of $\tilde\scrH$. 

The product $(j_1,j_2,j_3):=(h_1,h_2,h_3)\cdot(g_1,g_2,g_3)$ of two elements $(h_1,h_2,h_3),(g_1,g_2,g_3)$ of $\scrW_+(W(k))\times \scrW_0(W(k))\times\scrW_-(W(k))=\scrH(W(k))=\tilde\scrH(W(k))$, is defined by the identity $j_1j_2j_3^p=h_1h_2h_3^pg_1g_2g_3^p$.

The restrictions of $\scrP_{0-}$ to the factors $\scrW_+$, $\scrW_0$, and $\scrW_+\times_{W(k)} \scrW_0$ of $\scrH$ induce isomorphisms onto the closed subgroup schemes $\scrW_+$, $\scrW_0$, and $\scrW_{+0}$ (respectively) of $\pmb{GL}_M$. Moreover the restriction of $\scrP$ to the factors $\scrW_-$ and $\scrW_0\times_{W(k)} \scrW_-$ of $\scrH$ induces isomorphisms onto closed subgroup schemes of $\tilde\scrH$ which are isomorphic to $\scrW_-$ and $\scrW_{0-}$ (respectively). The product decomposition 
$$\tilde\scrH_{W_m(k)}=\scrH_{W_m(k)}=\scrW_{+W_m(k)}\times_{W_m(k)}\scrW_{0W_m(k)}\times_{W_m(k)} \scrW_{-W_m(k)}$$ 
into affine schemes depends on the reduction modulo $p^m$ of the chosen direct sum decomposition $M=F^1\oplus F^0$. But the group scheme $\tilde\scrH_{W_m(k)}=\scrH_{W_m(k)}$ over $\Spec(W_m(k))$ is intrinsically associated to $D$ i.e., it does not depend on the choice of the direct sum decomposition $M=F^1\oplus F^0$.

\medskip\noindent
{\bf 2.1.2. Coordinates.} To provide explicit descriptions of the affine schemes we have introduced, we will choose a $W(k)$-basis $\{e_1,\ldots,e_r\}$ for $M$ such that $\{e_1,\ldots,e_c\}$ is a $W(k)$-basis for $F^0$ and $\{e_{c+1},\ldots,e_r\}$ is a $W(k)$-basis for $F^1$. Then one can identify naturally $R_M=W(k)[x_{11},\ldots,x_{rr}][{1\over {\text{det}((x_{ij})_{1\le i,j\le r})}}]$, $R_{\tilde\scrH}=W(k)[p^{\eps_{11}}x_{11},\ldots,p^{\eps_{rr}}x_{rr}][{1\over {\text{det}((x_{ij})_{1\le i,j\le r})}}]$ where $\eps_{ij}\in\{-1,0\}$ is equal to $-1$ if and only if $1\le i\le c<j\le r$, $\scrW_{+}=\Spec(R_+)$ with $R_+:=W(k)[x_{ij}|1\le j\le c<i\le r]$, $\scrW_{0}=\Spec(R_0)$ with $R_0:=W(k)[x_{ij}|\text{either}\; 1\le i,j\le c\;\text{or}\; c<i,j\le r][{1\over {\text{det}((x_{ij})_{1\le i,j\le c})\text{det}((x_{ij})_{c< i,j\le r})}}]$, and $\scrW_{-}=\Spec(R_-)$ with $R_-:=W(k)[x_{ij}|1\le i\le c<j\le r]$.

\medskip\noindent
{\bf  2.1.3. Unipotent group schemes.} An affine group scheme over $k$ is called {\it unipotent} if it has a finite composition series whose factors are subgroup schemes of $\dbG_a$.  The class of unipotent group schemes over $k$ is stable under subgroup schemes, quotients,  and extensions (cf. [DG, Vol. II, Exp. XVII, Prop. 2.2]). A smooth, connected, affine group over $k$ is unipotent if and only it has no subgroup isomorphic to $\dbG_m$ (cf. [DG, Vol. II, Exp. XVII, Prop. 4.1.1]).

\medskip\noindent
{\bf  2.1.4. The $\dbW_m$ functor.} Let $\text{Aff}_k$ be the category of affine schemes over $k$. Let $\text{Set}$ and $\text{Group}$ be the categories of abstract sets and groups (respectively). Let $G$ be a smooth, affine scheme  of finite type (resp. a smooth, affine group scheme) over $\Spec(W(k))$. Let $\dbW_m(G):\text{Aff}_k\to\text{Set}$ (resp. $\dbW_m(G):\text{Aff}_k\to\text{Group}$) be the contravariant functor that associates to an affine $k$-scheme $\Spec(R)$ the set (resp. the group) $G(W_m(R))$. It is well known that this functor is representable by an affine scheme (resp. affine group scheme) over $k$ of finite type to be denoted also by $\dbW_m(G)$, cf. [Gr, Sect. 4, Cor. 1, p. 639] (resp. [Gr, Sect. 4, Cor. 4, p. 641]). We have $\dbW_m(G)(k)=G(W_m(k))$ and a natural identification $\dbW_1(G)=G_k$. If $I$ is an ideal of $R$ of square $0$, then the ideal $\Ker(W_m(R)\twoheadrightarrow W_m(R/I))$ is nilpotent and thus the reduction map $G(W_m(R))\to G(W_m(R/I))$ is surjective (cf. [BLR, Ch. 2, 2.2, Prop. 6]). From this and loc. cit. we get that the scheme (resp. the group scheme) $\dbW_m(G)$ is smooth. 

Suppose now that $G$ is a smooth, affine group scheme over $\Spec(W(k))$. If $m\ge 2$, then the length reduction epimorphisms $W_{m}(R)\twoheadrightarrow W_{m-1}(R)$ define naturally an epimorphism $\text{Red}_{m,G}:\dbW_{m}(G)\twoheadrightarrow \dbW_{m-1}(G)$ of groups over $k$. The kernel of $\text{Red}_{m,G}$ is the vector group over $k$ defined by $\Lie(G_k)$ and thus it is a unipotent, commutative group isomorphic to $\dbG_a^{\dim(G_k)}$. Using this and the identification $\dbW_1(G)=G_k$, by induction on $m\ge 1$ we get that (cf. also [Gr, Sect. 4, Cor. 4, p. 639] in connection to (iii)): 

\medskip
{\bf (i)} we have $\dim(\dbW_m(G))=m\dim(G_k)$;

\smallskip
{\bf (ii)} the group $\dbW_m(G)$ is connected if and only if $G_k$ is connected;

\smallskip
{\bf (iii)} if $H$ is a smooth, closed subgroup scheme of $G$, then the functorial homomorphism $\dbW_m(H)\to\dbW_m(G)$ is a closed embedding.

\medskip\smallskip\noindent
{\bf 2.2. The group action $\dbT_m$.} Let $\sigma_{\phi}:M\arrowsim M$ be the $\sigma$-linear automorphism which takes $x\in F^1$ to ${1\over p}\phi(x)$ and takes $x\in F^0$ to $\phi(x)$. Let $\sigma_{\phi}$ act on the sets underlying the groups $\pmb{GL}_M(W(k))$ and $\pmb{GL}_M(W_m(k))$ in the natural way: if $g\in\pmb{GL}_M(W(k))$, then $\sigma_{\phi}(g)=\sigma_{\phi} g\sigma_{\phi}^{-1}$ and $\sigma_{\phi}(g[m])=(\sigma_{\phi} g\sigma_{\phi}^{-1})[m]$. For $g\in\scrW_+(W(k))$ (resp. $g\in\scrW_0(W(k))$ or $g\in\scrW_-(W(k))$) we have $\phi(g)=\sigma_{\phi}(g^p)$ (resp. we have $\phi(g)=\sigma_{\phi}(g)$ or $\phi(g^p)=\sigma_{\phi}(g)$). 

Let $\scrH_m:=\dbW_m(\scrH)$ and $\scrD_m:=\dbW_m(\pmb{GL}_{M})$. 
As $\scrP_{W_m(k)}:\scrH_{W_m(k)}\to\tilde\scrH_{W_m(k)}$  
is an isomorphism of $\Spec(W_m(k))$-schemes, in all that follows we will identify naturally 
$$\scrH(W(k))=\tilde\scrH(W(k))\;\;\text{and}\;\;\scrH_m=\dbW_m(\scrH)=\dbW_m(\tilde\scrH).$$
Thus in what follows we will view $\scrH(W(k))$ as a subgroup of $\pmb{GL}_M(W(k))$ and $\scrH_m$ as a connected, smooth, affine group over $k$ of dimension $mr^2$ (cf. Subsubsection  2.1.4 applied to $\tilde\scrH$). Similarly, we will view $\scrD_m$ as a connected, smooth, affine variety over $k$ of dimension $mr^2$; occasionally (like in the proof of Theorem 2.4 (b) below), we will also view $\scrD_m$ as a smooth, affine group over $k$. We have a natural action 
$$\dbT_m:\scrH_m\times_k \scrD_m\to \scrD_m$$
defined on $k$-valued points as follows. If $h=(h_1,h_2,h_3)\in\scrH(W(k))$ and $g\in \pmb{GL}_{M}(W(k))$, then the product of $h[m]=(h_1[m],h_2[m],h_3[m])\in \scrH_m(k)=\scrH(W_m(k))$ and $g[m]\in\scrD_m(k)=\pmb{GL}_M(W_m(k))$ is the element
$$\dbT_m(h[m],g[m]):=(h_1h_2h_3^pg\phi(h_1h_2h_3^p)^{-1})[m]\leqno (1a)$$
$$=(h_1h_2h_3^pg\phi(h_3^p)^{-1}\phi(h_2)^{-1}\phi(h_1)^{-1})[m]=(h_1h_2h_3^pg\sigma_{\phi}(h_3)^{-1}\sigma_{\phi}(h_2)^{-1}\sigma_{\phi}(h_1^p)^{-1})[m]$$
$$=h_1[m]h_2[m]h_3[m]^pg[m]\sigma_{\phi}(h_3[m])^{-1}\sigma_{\phi}(h_2[m])^{-1}\sigma_{\phi}(h_1[m]^p)^{-1}\in\scrD_m(k).$$
The formula $\dbT_m(h[m],g[m])=(h_1h_2h_3^pg\phi(h_1h_2h_3^p)^{-1})[m]$ shows that the action $\dbT_m$ is intrinsically associated to $D$ i.e., it does not depend on the choice of the direct sum decomposition $M=F^1\oplus F^0$. For a later usage, we mention that
$$\dbT_1(h[1],g[1])=h_1[1]h_2[1]g[1]\sigma_{\phi}(h_3[1])^{-1}\sigma_{\phi}(h_2[1])^{-1}\in\scrD_1(k)=\pmb{GL}_M(k).\leqno (1b)$$ 
\indent
Let $\scrO_m$ be the orbit of $1_M[m]\in\scrD_m(k)$ under the action $\dbT_m$. Let $\bar\scrO_m$ be the Zariski closure of $\scrO_m$ in $\scrD_m$; it is an affine, integral scheme over $k$. The orbit $\scrO_m$ is a connected, smooth, open subscheme of $\bar\scrO_m$ and therefore it is also a quasi-affine scheme over $k$. Let $\scrS_m$ be the subgroup scheme of $\scrH_m$ which is the stabilizer of $1_M[m]$ under the action $\dbT_m$. Let $\scrC_m$ be the reduced group of $\scrS_m$. Let $\scrC_m^0$ be the identity component of $\scrC_m$. We have 
$$\dim(\scrS_m)=\dim(\scrC_m)=\dim(\scrC_m^0)=\dim(\scrH_m)-\dim(\scrO_m).\leqno (2)$$
\noindent
{\bf 2.2.1. Lemma.} {\it Let $g_1,g_2\in \pmb{GL}_M(W(k))$. Then the points $g_1[m],g_2[m]\in\scrD_m(k)$ belong to the same orbit of the action $\dbT_m$ if and only if the following two Dieudonn\'e modules $(M/p^mM,g_1[m]\phi_m,\vartheta_m g_1[m]^{-1})$ and  $(M/p^mM,g_2[m]\phi_m,\vartheta_m g_2[m]^{-1})$ are isomorphic.}

\medskip
\proof
Suppose that $g_1[m]$ and $g_2[m]$ belong to the same orbit of the action of $\dbT_m$. Let $h=(h_1,h_2,h_3)\in\scrH(W(k))$ be such that we have $\dbT_m(h[m],g_1[m])=g_2[m]$. Let $\break g_3:=h_1h_2h_3^pg_1\phi(h_1h_2h_3^p)^{-1}\in\pmb{GL}_M(W(k))$. The identity $\dbT_m(h[m],g_1[m])=g_2[m]$ implies that $g_3[m]=g_2[m]$. Thus $(M/p^mM,g_3[m]\phi_m,\vartheta_m g_3[m]^{-1})=(M/p^mM,g_2[m]\phi_m,\vartheta_m g_2[m]^{-1})$. As $(M,g_1\phi,\vartheta g_1^{-1})$ and $(M,g_3\phi,\vartheta g_3^{-1})$ are isomorphic, $(M/p^mM,g_1[m]\phi_m,\vartheta_m g_1[m]^{-1})$ and  $(M/p^mM,g_2[m]\phi_m,\vartheta_m g_2[m]^{-1})$ are also isomorphic. In particular, if $g_1[m]=g_2[m]$ we get that $(h_1h_2h_3^p)[m]$ is an automorphism of $(M/p^mM,g_1[m]\phi_m,\vartheta_m g_1[m]^{-1})$.

Suppose that $(M/p^mM,g_1[m]\phi_m,\vartheta_m g_1[m]^{-1})$ and $(M/p^mM,g_2[m]\phi_m,\vartheta_m g_2[m]^{-1})$ are isomorphic. Let $\tilde h\in \pmb{GL}_M(W(k))$ be such that $\tilde hg_1\phi(\tilde h)^{-1}$ and $g_2$ are elements of $\pmb{GL}_M(W(k))$ congruent  modulo $p^m$, cf. [Va2, Lem. 3.2.2] applied with $G=\pmb{GL}_M$. As $\tilde hg_1\phi(\tilde h)^{-1}\in \pmb{GL}_M(W(k))$, we have $\phi(\tilde h)\in \pmb{GL}_M(W(k))$. This implies that $\tilde h(\phi^{-1}(M))=\phi^{-1}(M)$. Thus $\tilde h(p\phi^{-1}(M))=p\phi^{-1}(M)$ and therefore $\tilde h(F^1+pM)=F^1+pM$. This implies that $\tilde h$  modulo $p$ normalizes $\bar F^1$ i.e., we have $\tilde h\in\scrH(W(k))=\tilde\scrH(W(k))$. Thus we can write $\tilde h=h_1h_2h_3^p$, where $h:=(h_1,h_2,h_3)\in\scrH(W(k))$. As $\tilde hg_1\phi(\tilde h)^{-1}$ and $g_2$ are congruent  modulo $p^m$, we have $\dbT_m(h[m],g_1[m])=g_2[m]$ i.e., $g_1[m]$ and $g_2[m]$ belong to the same orbit of the action $\dbT_m$.\endproof

\medskip
The following Corollary explains the title of this Section. 

\medskip\noindent
{\bf 2.2.2. Corollary.} {\it The set of orbits of the action $\dbT_m$ are in natural bijection to the set of isomorphism classes of truncated Barsotti--Tate groups of level $m$ over $k$ which have codimension $c$ and dimension $d$.}

\medskip
\proof
Let $\tilde B$ be a truncated Barsotti--Tate group of level $m$ over $k$ which has codimension $c$ and dimension $d$. Let $\tilde D$ be a $p$-divisible group over $k$ which lifts $\tilde B$, cf. [Il, Thm. 4.4 e)]. As $\tilde D$ has codimension $c$ and dimension $d$, its Dieudonn\'e module is isomorphic to $(M,\tilde g\phi)$ for some element $\tilde g\in\pmb{GL}_M(W(k))$. The Dieudonn\'e module of $\tilde B$ is isomorphic to $(M/p^mM,\tilde g[m]\phi_m,\vartheta_m\tilde g[m]^{-1})$. Based on this, the Corollary follows from Lemma 2.2.1 and the classical Dieudonn\'e theory (see [Fo, pp. 153 and 160]).\endproof

\medskip\noindent
{\bf 2.2.3. Lemma.} {\it We have a natural epimorphism $\Xi_m:\scrH_m\twoheadrightarrow\scrW_{0k}$ whose kernel is the unipotent radical $\scrH_m^{\text{unip}}$ of $\scrH_m$.}

\medskip
\proof The epimorphism $\Xi_m$ is defined at the level of $k$-valued points by the following rule: if $h[m]=(h_1[m],h_2[m],h_3[m])\in \scrH_m(k)$, then $\Xi_m(k)(h[m]):=h_2[1]\in\scrW_{0k}(k)$. We have a short exact sequence $1\to\Ker(\scrH_m\to\scrH_1)\to\Ker(\Xi_m)\to \scrW_{+k}\times_k\scrW_{-k}\to 1$, where the subgroup $\scrW_{+k}\times_k\scrW_{-k}$ of $\scrH_1$ is the usual product group.  As the group $\Ker(\scrH_m\to\scrH_1)$ has a composition series whose factors are the smooth, connected, unipotent groups $\Ker(\text{Red}_{l,\tilde\scrH})=\Ker(\scrH_l\to\scrH_{l-1})$ with $l\in\{2,\ldots,m\}$, we get that $\Ker(\Xi_m)$ has a composition series whose factors are smooth, connected, unipotent groups. Thus $\Ker(\Xi_m)$ is a smooth, connected, unipotent, normal subgroup of $\scrH_m$ and therefore it is a subgroup of  $\scrH_m^{\text{unip}}$. As $\scrW_{0k}=\pmb{GL}_{\bar F^1}\times_k\pmb{GL}_{\bar F^0}$ is a reductive group, $\scrH_m^{\text{unip}}$ is a subgroup of $\Ker(\Xi_m)$. We conclude that $\Ker(\Xi_m)=\scrH_m^{\text{unip}}$.\endproof

\medskip\smallskip\noindent
{\bf 2.3. On $\dbT_1$.} In this Subsection we assume that $m=1$ and we study the connected subgroup $\scrC_1^0$ of $\scrH_1=\scrH_k$ and $\dim(\scrO_1)$. In [Kr] (see also [Oo1, Subsect. (2.3) and Lem. (2.4)] and [Mo, Subsect. 2.1]) it is shown that there exists a $k$-basis $\{\bar e_1,\ldots,\bar e_r\}$ for $M/pM$ and a permutation $\pi$ of the set $J:=\{1,\ldots,r\}$ such that for $i\in J$ the following four properties hold: 

\medskip
{\bf (i)} $\phi_1(\bar e_i)=0$ if $i>c$; 

\smallskip
{\bf (ii)} $\phi_1(\bar e_i)=\bar e_{\pi(i)}$ if $i\le c$; 

\smallskip
{\bf (iii)} $\vartheta_1(\bar e_{\pi(i)})=0$ if $i\le c$; 

\smallskip
{\bf (iv)} $\vartheta_1(\bar e_{\pi(i)})=\bar e_i$ if $i>c$. 

\medskip
Let $\{e_1,\ldots,e_r\}$ be a $W(k)$-basis for $M$ that lifts the $k$-basis $\{\bar e_1,\ldots,\bar e_r\}$ for $M/pM$ and such that $F^1=\oplus_{i=c+1}^r W(k)e_i$. Let $\{e_{i,j}|i,j\in J\}$ be the $W(k)$-basis for $\End(M)$ such that for each $l\in J$ we have $e_{i,j}(e_l)=\delta_{j,l}e_i$. Let $\{\bar e_{i,j}|i,j\in J\}$ be the $k$-basis for $\End(M/pM)$ which is the reduction  modulo $p$ of $\{e_{i,j}|i,j\in J\}$. Let $\sigma_{\pi}:M\arrowsim M$ be the $\sigma$-linear automorphism that takes $e_i$ to $e_{\pi(i)}$ for all $i\in J$. Let $g_{\pi}:=\sigma_{\pi}\sigma_{\phi}^{-1}\in\pmb{GL}_M(W(k))$. Due to the properties (i) to (iv), the reduction modulo $p$ of $(M,g_{\pi}\phi,\vartheta g_{\pi}^{-1})$ coincides with $(M/pM,\phi_1,\vartheta_1)$. Based on this and Lemma 2.2.1, (up to isomorphisms) we can assume that $g_{\pi}[1]=1_M[1]$; thus $\sigma_{\phi}$ and $\sigma_{\pi}$ are congruent modulo $p$. As the action $\dbT_1$ is intrinsically associated to $D$ (i.e., it does not depend on the choice of the direct sum decomposition $M=F^1\oplus F^0$), to study the group $\scrC_1^0$ we can assume that $F^0=\oplus_{i=1}^c W(k)e_i$. We define 
$$\scrJ_+:=\{(i,j)\in J^2|j\le c<i\},\;\; $$ 
$$\scrJ_0:=\{(i,j)\in J^2|\;\text{either}\;i,j>c\;\text{or}\;i,j\le c\},\;\;\;\text{and}\;\;\;\scrJ_-:=\{(i,j)\in J^2|i\le c<j\}.$$
\noindent
The three sets $\{\bar e_{i,j}|(i,j)\in\scrJ_+\}$, $\{\bar e_{i,j}|(i,j)\in\scrJ_0\}$, and $\{\bar e_{i,j}|(i,j)\in\scrJ_-\}$ are $k$-basis for $\Lie(\scrW_{+k})=\Hom(\bar F^0,\bar F^1)$, $\Lie(\scrW_{0k})=\End(\bar F^1)\oplus\End(\bar F^0)$, and $\Lie(\scrW_{-k})=\Hom(\bar F^1,\bar F^0)$ (respectively).   

Let $(h_1[1],h_2[1],h_3[1])\in\scrH_1(k)$. We have $(h_1[1],h_2[1],h_3[1])\in\scrC_1(k)$ if and only if the following identity holds (cf. (1b))
$$h_1[1]h_2[1]=\sigma_{\phi}(h_2[1])\sigma_{\phi}(h_3[1]).\leqno(3)$$ 
Let $h_{12}[1]:=h_1[1]h_2[1]\in\scrW_{+0}(k)$ and $h_{23}[1]:=h_2[1]h_3[1]\in\scrW_{0-}(k)$. We write 
$$h_{12}[1]=1_M[1]+\sum_{(i,j)\in\scrJ_+\cup\scrJ_0} x_{i,j}\bar e_{i,j\;}\;\text{and}\;\;h_{23}[1]=1_M[1]+\sum_{(i,j)\in\scrJ_0\cup\scrJ_-} x_{i,j}\bar e_{i,j}.$$ 
We have $\sigma_{\pi}(x_{i,j}\bar e_{i,j})=x_{i,j}^p\bar e_{\pi(i),\pi(j)}$. Thus, as $\sigma_{\phi}$ and $\sigma_{\pi}$ are congruent modulo $p$, we get that the identity (3) (i.e., the identity $h_{12}[1]=\sigma_{\phi}(h_{23}[1])$) holds if and only if for all pairs $(i,j)\in J^2$ we have:

\medskip
{\bf (v)}  $x_{\pi(i),\pi(j)}=x_{i,j}^p$, if $(i,j)\in\scrJ_-\cup\scrJ_0$ and $(\pi(i),\pi(j))\in\scrJ_+\cup\scrJ_0$;

\smallskip
{\bf (vi)} $x_{\pi(i),\pi(j)}=0$, if $(i,j)\in\scrJ_+$ and $(\pi(i),\pi(j))\in\scrJ_+\cup\scrJ_0$;

\smallskip
{\bf (vii)} $0=x_{i,j}^p$, if $(i,j)\in\scrJ_-\cup\scrJ_0$ and $(\pi(i),\pi(j))\in\scrJ_-$.

\medskip
Let $\scrJ_-^{\pi}$ be the subset of $\scrJ_-$ formed by those pairs $(i,j)$ with the property that if $\nu_{\pi}(i,j)$ is the smallest positive integer such that 
$$(\pi^{\nu_{\pi}(i,j)}(i),\pi^{\nu_{\pi}(i,j)}(j))\in\scrJ_+\cup\scrJ_-,$$ 
then we have $(\pi^{\nu_{\pi}(i,j)}(i),\pi^{\nu_{\pi}(i,j)}(j))\in\scrJ_+$. Based on the properties (v) to (vii) we easily get that the variable $x_{\tilde i,\tilde j}$ can take an infinite number of values if and only if there exist a pair $(i,j)\in\scrJ_-^{\pi}$ and a number $l\in\{0,1,\ldots,\nu_{\pi}(i,j)\}$ such that we have $(\tilde i,\tilde j)=(\pi^l(i),\pi^l(j))$. This implies that we have $(h_1[1],h_2[1],h_3[1])\in\scrC_1^0(k)$ if and only if we have identities
$$h_{12}[1]=1_M[1]+\sum_{(i,j)\in\scrJ_-^{\pi}} \sum_{l=1}^{\nu_{\pi}(i,j)} x_{i,j}^{p^l}\bar e_{\pi^l(i),\pi^l(j)},\leqno (4a)$$
$$h_2[1]=1_M[1]+\sum_{(i,j)\in\scrJ_-^{\pi}}\sum_{l=1}^{\nu_{\pi}(i,j)-1} x_{i,j}^{p^l}\bar e_{\pi^l(i),\pi^l(j)},\leqno (4b)$$
$$h_{23}[1]=1_M[1]+\sum_{(i,j)\in\scrJ_-^{\pi}} \sum_{l=0}^{\nu_{\pi}(i,j)-1} x_{i,j}^{p^l}\bar e_{\pi^l(i),\pi^l(j)},\leqno (4c)$$where the variables $x_{i,j}$ with $(i,j)\in\scrJ_-^{\pi}$ can take independently all those values in $k$ for which we have $h_2[1]\in\scrW_0(k)$. Based on Formulas (4a) to (4c) we get that 
$$\Lie(\scrC_1^0)=\oplus_{(i,j)\in\scrJ_-^{\pi}} k\bar e_{i,j}\subseteq\Lie(\scrW_{-k}).\leqno (4d)$$
\indent
From Formula (4d) we get that each element of $\Lie(\scrC_1^0)$ is a nilpotent element of $\End(\bar M)$. Thus $\Lie(\scrC_1^0)$ has no non-zero semisimple element. Therefore $\scrC_1^0$ has no subgroup isomorphic to $\dbG_m$ and thus (cf. Subsubsection 2.1.4) we have:

\medskip
{\bf (viii)} the smooth, connected, affine group $\scrC_1^0$ is unipotent.

\medskip
From Formula (4d) we also get that:

\medskip
{\bf (ix)} the dimension of $\scrC_1^0$ is the number of elements of the subset $\scrJ_-^{\pi}$ of $\scrJ_-$.

\medskip
As $\dim(\scrH_1)=r^2$ is equal to the number of elements of $J^2$, from the property (ix) and Formula (2) applied with $m=1$ we get that:

\medskip
{\bf (x)} the dimension of $\scrO_1$ is  the number of elements of the set $\scrJ^2\setminus\scrJ_-^{\pi}$.

\medskip

The next two examples show how the property (ix) can be used to compute $\dim(\scrC_1^0)$. 

\medskip\noindent
{\bf 2.3.1. Example.} In this example all indices $i,j\in J=\{1,\ldots,r\}$ are taken  modulo $r$. Suppose that $v=1$ (thus $g.c.d.(c,d)=1$) and that $D$ is {\it minimal} in the sense of [Oo3, Subsect. (1.1)]. We can assume that we have $\pi(i)=i+d$ for all $i\in J$, cf. [Oo3, Subsect. (1.4)]. The $p$-divisible group $D$ is uniquely determined up to isomorphism by $D[p]$, cf. [Va2, Example 3.3.6]. Thus $n_D=1$ and therefore the Dieudonn\'e modules $(M,\phi)$ and $(M,g_{\pi}\phi)$ are isomorphic. Thus we can assume that $g_{\pi}=1_M$ i.e., $\sigma_{\phi}=\sigma_{\pi}$. Therefore for all $i,j\in J$ we have $\phi(e_{i,j})=p^{n_{i,j}}e_{i+d,j+d}$, where the number $n_{i,j}\in\{-1,0,1\}$ is defined by the rule:

\medskip
{\bf (*)} it is $-1$ if $(i,j)\in \scrJ_-$, it is $0$ if $(i,j)\in \scrJ_0$, and it is $1$ if $(i,j)\in \scrJ_+$.

\medskip
Suppose that $(i,j)\in\scrJ_-$. From [Va2, Example 3.3.6] we get that the first non-zero element of the sequence $(n_{i+ld,j+ld})_{l\ge 1}$ is $1$; thus $(i,j)\in\scrJ_-^{\pi}$. In other words, we have $\scrJ_-^{\pi}=\scrJ_-$. As $\scrJ_-$ has $cd$ elements, we get that $\dim(\scrC_1^0)=cd$ (cf. property 2.3 (ix)).  

\medskip\noindent
{\bf 2.3.2. Example.} Suppose that $v=2$ and that $\alpha_1<\alpha_2$. Thus $c=c_1+c_2$, $d=d_1+d_2$, and $c_2d_1<c_1d_2$. Suppose that we have a product decomposition $D=D_1\times_k D_2$ such that the heights of $D_1$ and $D_2$ are $r_1$ and $r_2$ (respectively) and both $D_1$ and $D_2$ are minimal $p$-divisible groups over $k$. The $p$-divisible group $D$ is uniquely determined up to isomorphism by $D[p]$, cf. either [Oo3, Thm. 1.2] or [Va4, Main Thm. D]. Thus $n_D=1$ and therefore as in Example 2.3.1 we argue that we can assume that $\sigma_{\phi}=\sigma_{\pi}$. 

In what follows, the indices $i$ of the letter $a$ are taken modulo $r_2$ and the indices $j$ of the letter $b$ are taken modulo $r_1$. Let $(a_1,\ldots,a_{r_2}):=(e_{c_1}+1,\ldots,e_{c},e_{r_1+c_2+1},\ldots,e_r)$ and $(b_1,\ldots,b_{r_1}):=(e_1,\ldots,e_{c_1},e_{c+1},\ldots,e_{r_1+c_2})$. We can assume that the permutation $\pi$ of $J$ is such that we have $\sigma_{\pi}(a_i)=a_{i+d_2}$ for all $i\in\{1,\ldots,r_2\}$ and we have $\sigma_{\pi}(b_j)=b_{j+d_1}$ for all $j\in\{1,\ldots,r_1\}$, cf. Example 2.3.1 applied to $D_2$ and $D_1$ (respectively). Let $M_1:=\oplus_{j=1}^{r_1} W(k)b_j$ and $M_2:=\oplus_{i=1}^{r_2} W(k)a_i$. The pairs $(M_1,\phi)$ and $(M_2,\phi)$ are the Dieudonn\'e modules of $D_1$ and $D_2$ (respectively). We consider the disjoint union decomposition 
$$\scrJ_-=\scrJ_{-,1,1}\cup\scrJ_{-,1,2}\cup\scrJ_{-,2,1}\cup\scrJ_{-,2,2}\leqno (5a)$$ 
such that the pair $(i,j)\in\scrJ_-$ belongs to $\scrJ_{-,1,1}$ (resp. to $\scrJ_{-,1,2}$, $\scrJ_{-,2,1}$, or $\scrJ_{-,2,2}$) if and only if $e_{i,j}$ belongs to $\End(M_1)$ (resp. to $\Hom(M_1,M_2)$, $\Hom(M_2,M_1)$, or $\End(M_2)$). To the decomposition (5) corresponds a disjoint union decomposition 
$$\scrJ_-^{\pi}=\scrJ_{-,1,1}^{\pi}\cup\scrJ_{-,1,2}^{\pi}\cup\scrJ_{-,2,1}^{\pi}\cup\scrJ_{-,2,2}^{\pi},\leqno (6)$$
 where $\scrJ_{-,\dag,\ddag}^{\pi}:=\scrJ_{-,\dag,\ddag}\cap \scrJ_-^{\pi}$. From Example 2.3.1 applied to $D_1$ and $D_2$, we get that the sets $\scrJ^{\pi}_{-,1,1}$ and $\scrJ^{\pi}_{-,2,2}$ have $c_1d_1$ and $c_2d_2$ (respectively) elements. 

We check that each one of the two sets $\scrJ_{-,1,2}^{\pi}$ and $\scrJ_{-,2,1}^{\pi}$ has $c_2d_1$ elements. We will perform this computation only for $\scrJ_{-,1,2}^{\pi}$ as the computation for $\scrJ_{-,2,1}^{\pi}$ is the same (via standard duality in which the roles of $-$ and $+$ are interchanged). If $\{e_1^*,\ldots,e_r^*\}$ is the $W(k)$-basis for $\Hom(M,W(k))$ that is the dual of the $W(k)$-basis $\{e_1,\ldots,e_r\}$ for $M$, then under the identification $\End(M)=M\otimes_{W(k)} \Hom(M,W(k))$ we have $e_{i,j}=e_i\otimes e_j^*$ for all $i,j\in J$. As $\sigma_{\phi}=\sigma_{\pi}$, for $(i,j)\in \{1,\ldots,r_2\}\times\{1,\ldots,r_1\}$ we have $\phi(a_i\otimes b_j^*)=p^{n_{i,j}}a_{i+d_2}\otimes b_{j+d_1}^*$, where the number $n_{i,j}\in\{-1,0,1\}$ is defined by the rule:

\medskip
{\bf (*)} it is $-1$ if $i\le c_2$ and $j>c_1$ (i.e., if $(i,j)\in\scrJ_{-,1,2}$), it is $0$ if $i\le c_2$ and $j\le c_1$ or if $i>c_2$ and $j>c_1$, and it is $1$ if $i>c_2$ and $j\le c_1$.
\medskip\noindent
The number of elements of the set $\scrJ_{-,1,2}^{\pi}$ is the number of pairs $(i,j)\in \{1,\ldots,r_2\}\times\{1,\ldots,r_1\}$ with the properties that $n_{i,j}=-1$ and that the first non-zero element of the sequence $(n_{i+ld_2,j+ld_1})_{l\ge 1}$ is $1$. The set $\{(i,j)\in \{1,\ldots,r_2\}\times\{1,\ldots,r_1\}|n_{i,j}=-1\}$ has $c_2d_1$ elements. Thus to prove that the set $\scrJ_{-,1,2}^{\pi}$ has $c_2d_1$ elements, it suffices to show that there exists no pair $(i,j)\in \{1,\ldots,r_2\}\times\{1,\ldots,r_1\}$ such that  $n_{i,j}=-1$ and the first non-zero element of the sequence $(n_{i+ld_2,j+ld_1})_{l\ge 1}$ is $-1$. This is an elementary number theory property that can be checked using the ideas of [Va2, Example 3.3.6]. For a change, here we will check this property using the ideas of [Va4]. 

We show that the assumption that there exists a pair $(i,j)\in \{1,\ldots,r_2\}\times\{1,\ldots,r_1\}$ such that  $n_{i,j}=-1$ and the first non-zero element of the sequence $(n_{i+ld_2,j+ld_1})_{l\ge 1}$ is $-1$, leads to a contradiction. Let $q$ be the positive integer such that $n_{i,j}=n_{i+qd_2,j+qd_1}=-1$ and $n_{i+d_2,j+d_1}=\cdots=n_{i+(q-1)d_2,j+(q-1)d_1}=0$. The element
$$g_{i,j}:=1_M+pa_i\otimes b_j^*\in\Ker(\pmb{GL}_M(W(k))\to \pmb{GL}_M(k))$$ 
fixes all elements of $\{e_1,\ldots,e_r\}\setminus\{b_j\}$ and takes $b_j$ to $b_j+pa_i$. As $n_D=1$, the Dieudonn\'e modules $(M,\phi)$ and $(M,g_{i,j}\phi)$ are isomorphic. Thus there exists an element $h\in \pmb{GL}_M(W(k))$ such that $g_{i,j}\phi=h\phi h^{-1}$. As $g_{i,j}$ acts trivially on $M_2$ and $M/M_2$, the restriction of $h$ to $M_2$ is an automorphism $h_2$ of $(M_2,\phi)$ and the $W(k)$-automorphism of $M/M_2$ induced by $h$ is an automorphism $h_1$ of $(M/M_2,\phi)$. We identify $h_1$ with an automorphism of $(M_1,\phi)$. By replacing $h$ with $h(h_1\oplus h_2)^{-1}$, we can assume that $h$ acts trivially on $M_2$ and $M/M_2$. Therefore we can write $h=1_M+u$, where $u\in\Hom(M_1,M_2)$. 

As $g_{i,j}\phi=h\phi h^{-1}$, we have $g_{i,j}=h\phi(h^{-1})=(1_M+u)(1_M-\phi(u))=1_M+u-\phi(u)$; thus $pa_i\otimes b_j^*=u-\phi(u)$. As $\alpha_1<\alpha_2$, all Newton polygon slopes of $(\Hom(M_1,M_2)[{1\over p}],\phi)$ are positive. Thus the equation $pa_i\otimes b_j^*=u-\phi(u)$ has a unique solution $u=\sum_{l=0}^{\infty} \phi^l(pa_i\otimes b_j^*)$. We write $u=\sum_{s=1}^{r_2}\sum_{t=1}^{r_1} u_{s,t}a_s\otimes b_t^*$, where each $u_{s,t}\in W(k)$. As all Newton polygon slopes of $(\Hom(M_1,M_2)[{1\over p}],\phi)$ are positive, the pairs $(i+ld_2,j+ld_1)$ with $l\in\{0,\ldots,q\}$ are distinct and the smallest positive integer $q_0$ such that $(i,j)=(i+q_0d_2,j+q_0d_1)$ is at least equal to $q+4$. Moreover (as $\sigma_{\phi}=\sigma_{\pi}$), for each pair $(s,t)\in\{(i+ld_2,j+ld_1)|l\in\{0,\ldots,q_0-1\}\}$ we have $\phi^{q_0}(a_s\otimes b_t^*)=p^{m_{s,t}}a_s\otimes b_t^*$ for some positive integer $m_{s,t}$. As $\sum_{l=0}^{q} n_{i+ld_2,j+ld_1}=-2$, we have $\phi^{q+1}(pa_i\otimes b_j^*)={1\over p}a_{i+(q+1)d_2}\otimes b_{j+(q+1)d_1}^*$. As $u=\sum_{l=0}^{\infty} \phi^l(pa_i\otimes b_j^*)$, from the last two sentences we get that 
$$u_{i+(q+1)d_2,j+(q+1)d_1}={1\over p}\sum_{l=0}^{\infty} p^{lm_{i+(q+1)d_2,j+(q+1)d_1}}\notin W(k).$$
Contradiction. 

Thus each one of the two sets $\scrJ_{-,1,2}^{\pi}$ and $\scrJ_{-,2,1}^{\pi}$ has $c_2d_1$ elements.
Due to (6), we conclude that the set $\scrJ_-^{\pi}$ has $c_1d_1+c_2d_2+2c_2d_1$ elements. From the property 2.3 (ix) we get that $\dim(\scrC_1^0)=c_1d_1+c_2d_2+2c_2d_1$.

\medskip\smallskip\noindent
{\bf 2.4. Theorem.} {\it  We recall that $m$ is a positive integer. Let $\pmb{\text{Aut}}(D[p^m])_{\text{crys},\red}$ be the reduced group of the group scheme $\pmb{\text{Aut}}(D[p^m])_{\text{crys}}$ over $k$ of automorphisms of the triple $(M/p^mM,\phi_m,\vartheta_m)$. Then the following three properties hold:

\medskip
{\bf (a)} the connected, smooth group $\scrC_m^0$ is unipotent;

\smallskip
{\bf (b)} there exist two finite epimorphisms 
$$\iota_m:\scrC_m\twoheadrightarrow\pmb{\text{Aut}}(D[p^m])_{\text{crys},\red}\;\;\text{and}\;\;\lambda_{m,\red}:\pmb{\text{Aut}}(D[p^m])_{\red}\twoheadrightarrow \pmb{\text{Aut}}(D[p^m])_{\text{crys},\red}$$ which at the level of $k$-valued points induce isomorphisms $\iota_m(k):\scrC_m(k)\arrowsim\pmb{\text{Aut}}(D[p^m])_{\text{crys},\red}(k)$ and $\lambda_{m,\red}(k):\pmb{\text{Aut}}(D[p^m])_{\red}(k)\arrowsim \pmb{\text{Aut}}(D[p^m])_{\text{crys},\red}(k)$;

\smallskip
{\bf (c)} we have $\dim(\scrS_m)=\dim(\scrC_m)=\dim(\scrC^0_m)=\gamma_D(m)$.}

\medskip
\proof
We prove (a) by induction on $m\ge 1$. The basis of the induction holds for $m=1$, cf. property 2.3 (viii). For $m\ge 2$, the passage from $m-1$ to $m$ goes as follows. The epimorphism $\text{Red}_{m,\tilde\scrH}:\scrH_m\twoheadrightarrow\scrH_{m-1}$ maps the group $\scrC_m^0$ to $\scrC_{m-1}^0$. The group $\Ker(\text{Red}_{m,\tilde\scrH})$ is unipotent, cf. Subsubsection  2.1.4. The group schemes $\im(\scrC_m^0\to\scrC_{m-1}^0)\leqslant \scrC_{m-1}^0$ and $\Ker(\scrC_m^0\to \scrC_{m-1})\leqslant \Ker(\text{Red}_{m,\tilde\scrH})$ are also unipotent, cf. Subsubsection  2.1.3. As $\scrC_m^0$  is the extension of $\im(\scrC_m^0\to \scrC_{m-1}^0)$ by $\Ker(\scrC_m^0\to\scrC_{m-1}^0)$, we get that $\scrC_m^0$ is a unipotent group (cf. Subsubsection  2.1.3). This ends the induction. Thus (a) holds.

To prove (b), we view $\scrD_m$ as a connected, smooth, affine group over $k$. If $R$ is a commutative $k$-algebra, then $\pmb{\text{Aut}}(D[p^m])_{\text{crys}}(R)$ is the subgroup of $\scrD_m(R)=\pmb{GL}_M(W_m(R))$ formed by elements that commute with $\phi_m\otimes 1_{W_m(R)}$ and $\vartheta_m\otimes 1_{W_m(R)}$; here $\phi_m$ and $\vartheta_m$ are viewed as $W_m(k)$-linear maps $M/p^mM\otimes_{W_m(k)}{}_{\sigma} W_m(k)\to M/p^mM$ and $M/p^mM\to M/p^mM\otimes_{W_m(k)}{}_{\sigma} W_m(k)$ (respectively). Thus $\pmb{\text{Aut}}(D[p^m])_{\text{crys}}$ is a subgroup scheme of $\scrD_m$. The crystalline Dieudonn\'e theory provides us with a homomorphism 
$$\lambda_m:\pmb{\text{Aut}}(D[p^m])\to \pmb{\text{Aut}}(D[p^m])_{\text{crys}}$$
that takes an element $x\in \pmb{\text{Aut}}(D[p^m])(R)$ to the inverse of the element of $\pmb{\text{Aut}}(D[p^m])_{\text{crys}}(R)$ which is the evaluation of $\dbD(x)$ at the thickening $\grS_m(R)$ (here we need the inverse, as we use contravariant Dieudonn\'e modules and thus theory). 

Let $\lambda_{m,\red}:\pmb{\text{Aut}}(D[p^m])_{\text{red}}\to \pmb{\text{Aut}}(D[p^m])_{\text{crys},\red}$
be the homomorphism between reduced groups defined by $\lambda_m$. The  homomorphism $\lambda_m(k):\pmb{\text{Aut}}(D[p^m])(k)\to \pmb{\text{Aut}}(D[p^m])_{\text{crys}}(k)$ is an isomorphism, cf. classical Dieudonn\'e theory. We can naturally identify $\lambda_{m,\red}(k)$ with $\lambda_m(k)$. Thus the homomorphism $\lambda_{m,\red}(k):\pmb{\text{Aut}}(D[p^m])_{\red}(k)\to \pmb{\text{Aut}}(D[p^m])_{\text{crys},\red}(k)$ is an isomorphism. This implies that the homomorphism $\lambda_{m,\red}$ is a finite epimorphism.

Let $\iota_m:\scrC_m\to \pmb{\text{Aut}}(D[p^m])_{\text{crys},\red}$ be the homomorphism which takes an element $h[m]=(h_1[m],h_2[m],h_3[m])\in\scrC_m(k)\leqslant\scrH_m(k)$ to the element $h_1[m]h_2[m]h_3^p[m]\in\scrD_m(k)=\pmb{GL}_M(W_m(k))$. The fact that $h_1[m]h_2[m]h_3[m]^p\in \pmb{\text{Aut}}(D[p^m])_{\text{crys},\red}(k)$ was already checked at the end of the first paragraph of the proof of Lemma 2.2.1. If we have $h[m]\in \Ker(\iota_m(k))$, then $h_1[m]=1_M[m]$, $h_2=1_M[m]$, and $h_3[m]^p=1_M[m]$. As $h[m]\in\scrC_m(k)$, the element $\dbT_m(h[m],1_M[m])=1_M[m]\sigma_{\phi}(h_3[m])^{-1}$ is $1_M[m]$; therefore $\sigma_{\phi}(h_3[m])=1_M[m]$. Thus $h_3[m]=1_M[m]$. Therefore $h[m]\in\Ker(\iota_m(k))$ is the identity element. Thus the kernel of $\iota_m(k)$ is trivial. 

Next we check that $\iota_m(k)$ is onto. Let $\tilde h\in\pmb{\text{Aut}}(D[p^m])_{\text{crys},\red}(k)\leqslant\pmb{GL}_M(W_m(k))$. As $\tilde h$  modulo $p$ normalizes $\bar F_1$, we can write $\tilde h=h_1[m]h_2[m]h_4[m]^p$ for some element $(h_1,h_2,h_4)\in\scrH(W(k))$. Let $g:=h_1h_2h_4^p\phi(h_1h_2h_4^p)^{-1}\in\pmb{GL}_M(W(k))$. As we have $\tilde h\in\pmb{\text{Aut}}(D[p^m])_{\text{crys},\red}(k)$, we can identify $(M/p^mM,g[m]\phi_m,\vartheta_m g[m]^{-1})=(M/p^mM,\phi_m,\vartheta_m)$. As $g[m]\phi_m=\phi_m$, the element $g[m]$ fixes $\phi(M)/p^mM$. Thus we can write $g=1_M+p^{m-1}e$, where $e\in\End(M)$ is such that $e$  modulo $p$ annihilates $\phi(M)/pM$. As $\vartheta_m=\vartheta_mg[m]^{-1}$, the reduction modulo $p$ of $e$ annihilates $\Ker(\vartheta_1)=M/\phi(M)$. As $e$ modulo $p$ annihilates both $\phi(M)/M=\sigma_{\phi}(F^1)+pM/pM$ and $M/\phi(M)=M/\sigma_{\phi}(F^1)+pM$, the reduction  modulo $p$ of $1_M+e$ belongs to $(\sigma_{\phi}\scrW_-\sigma_{\phi}^{-1})(k)$. Thus there exists an element $h_5\in\Ker(\scrW_-(W_m(k))\to\scrW_-(W_{m-1}(k)))$ such that we have $g[m]=\sigma_{\phi}(h_5)[m]$. Let $h:=(1_M,1_M,h_5)\cdot (h_1,h_2,h_4)\in\scrH(W(k))$. As $h_5[m]^p=1_M[m]$, we have $h[m]\in\scrC_m(k)$ due to the following identities
$$\dbT_m(h[m],1_M[m])=(h_1h_2h_4^p\phi(h_1h_2h_4^p)^{-1}\sigma_{\phi}(h_5)^{-1})[m]=(h_1h_2h_4^p\phi(h_1h_2h_4^p)^{-1}g^{-1})[m]=1_M[m].$$ 
As $h_5[m]^p=1_M[m]$, we have $\iota_m(h[m])=(h_5^ph_1h_2h_4^p)[m]=(h_1h_2h_4^p)[m]=\tilde h$. Thus $\iota_m(k)$ is onto and therefore $\iota_m(k)$ is an isomorphism. This implies that $\iota_m$ is a finite epimorphism. Thus (b) holds. 

Due to Formula (2), (c) follows from (b) (see Definition 1.1 (a) for $\gamma_D(m)$).\endproof

\medskip\noindent
{\bf 2.4.1. Remark.} In [Tr2, Thm. 2] and [Tr3, Sect. 26], Traverso considered a group action that is very close in nature to the action $\dbT_m$ and that has the form 
$$\dbT_m^{\text{old}}:\scrD_m\times_k\scrV_m\to\scrV_m,$$ where $\scrV_m$ is an affine variety over $k$ of dimension $rm^2-d^2$. We emphasize that:

\medskip\noindent
-- the description of $\scrV_m$ is less clear and usable than the descriptions of $\scrD_m$ and $\scrH_m$ (see [Tr2, Subsect. 2.5] and [Tr3, Sect. 26]);

\smallskip\noindent
-- the stabilizer subgroup schemes of the action $\dbT_m^{\text{old}}$ do not have a nice geometric interpretation (to be compared with Theorem 2.4 (b) and (c)); more precisely, the $k$-valued groups of the resulting stabilizer subgroup schemes are of the form $\{g[m]\in\scrD_m(k)|g[m]\phi_m=\phi_m g[m]\}$ (cf. [Tr2, proof of Subsect. 2.6, p. 51]) and therefore they ignore the Verschiebung maps (like $\vartheta_m$);

\smallskip\noindent
-- the action $\dbT_m^{\text{old}}$ does not generalize easily to relative contexts (like the ones of Shimura $F$-crystals used in [Va3] or the more general ones to be introduced in Subsection 4.1 below) or to the context of $p$-divisible objects over $k$ (see [Va2, Subsubsect. 2.2.1 (d)]).  

\medskip\smallskip\noindent
{\bf 2.5. Lemma.} {\it Suppose that the image of the abstract composite homomorphism 
$$\chi_D(m):\pmb{\text{Aut}}(D[p^m])(k)\to\pmb{\text{Aut}}(D[p])(k)\to\scrW_{+0}(k)\twoheadrightarrow (\scrW_{+0}/\scrW_{+})(k)\arrowsim\scrW_0(k)$$ 
is finite (the homomorphism $\pmb{\text{Aut}}(D[p])(k)\to\scrW_{+0}(k)$ is defined naturally by the homomorphism $\lambda_{1,\text{red}}$ of Theorem 2.4 (b)). Then the connected group $\scrC_m^0$ is a subgroup of $\scrH_m^{\text{unip}}$.}

\medskip
\proof 
The homomorphisms $\iota_m$ and $\lambda_{m,\text{red}}$ of Theorem 2.4 (b) are compatible with the standard reduction homomorphisms $\scrC_m\to\scrC_1$, $\pmb{\text{Aut}}(D[p^m])_{\text{crys},\red}\to\pmb{\text{Aut}}(D[p])_{\text{crys},\red}$, and $\pmb{\text{Aut}}(D[p^m])_{\red}\to\pmb{\text{Aut}}(D[p])_{\red}$. Thus from Theorem 2.4 (b) and our hypotheses we get that the homomorphism $\scrC_m\to\scrW_{0k}$ which is the restriction of the homomorphism $\Xi_m$ of Lemma 2.2.3, has finite image. Therefore the connected group $\scrC_m^0$ is a subgroup of $\Ker(\Xi_m)$. From this and Lemma 2.2.3 we get that $\scrC_m^0$ is a subgroup of $\scrH_m^{\text{unip}}$. 
\endproof

\medskip\smallskip\noindent
{\bf 2.6. Lemma.} {\it We assume that $\scrC_m^0$ is a subgroup of $\scrH_m^{\text{unip}}$. Then the orbit $\scrO_m$ (of $1_M[m]$ under the action $\dbT_m$) is an affine scheme.}

\medskip
\proof
The natural epimorphism $\scrH_m/\scrC_m^0\twoheadrightarrow\scrH_m/\scrS_m\arrowsim\scrO_m$ is finite. Thus from Chevalley's Theorem (see [Gro, Ch. II, (6.7.1)]) we get that $\scrO_m$ is affine if and only if $\scrH_m/\scrC_m^0$ is affine. The homogeneous space  $\scrH_m/\scrC_m^0$ is affine if and only if $\scrC_m^0$ is an {\it exact} subgroup of $\scrH_m$, cf. [CPS, Thm. 4.3]. We recall (cf. [CPS]) that a subgroup $\triangle$ of $\scrH_m$ is called exact if and only if the induction of rational $\triangle$-modules to rational $\scrH_m$-modules preserves short exact sequences. As $\scrC_m^0$ is a subgroup of $\scrH_m^{\text{unip}}$, it is an exact subgroup of $\scrH_m$ (cf. [CPS, Cor. 4.6]). Therefore the homogeneous space $\scrH_m/\scrC_m^0$ is affine. Thus $\scrO_m$ is affine.\endproof

\bigskip
\noindent
{\boldsectionfont 3. The proof of the Basic Theorem}
\bigskip 
 
In this Section we prove the Basic Theorem. We use the notations listed before Subsection 2.1. The notations $M=F^1\oplus F^0$, $\bar F^1$, $\bar F^0$, $\scrW_+$, $\scrW_0$, $\scrW_-$, $\scrW_{+0}$, $\scrW_{0-}$, $\scrH$, $\tilde\scrH$, $\phi_m,\vartheta_m:M/p^mM\to M/p^mM$, $\scrP_{0-W_m(k)}=\tilde\scrP_{0-W_m(k)}$, $\sigma_{\phi}$, $\scrH_m$, $\scrD_m$, $\dbT_m$, $\scrO_m$, $\scrC_m$, and $\scrC_m^0$ are as in Subsections 2.1 and 2.2. If $\bullet$ is a commutative, flat $W(k)$-algebra, let $\Omega_{\bullet}^\wedge$ be the $p$-adic completion of the $\bullet$-module $\Omega_{\bullet}$ of differentials of $\bullet$ and let $\delta_0:M\otimes_{W(k)} \bullet\to M\otimes_{W(k)} \Omega_{\bullet}^\wedge$ be the flat connection that annihilates $M\otimes 1$.

\medskip\smallskip\noindent
{\bf 3.1. Proof of 1.2 (a).}
To prove Theorem 1.2 (a), we can work locally in the Zariski topology of $\scrA$. Thus fixing a point $y_0\in\scrA(k)$, we can assume that $\scrA=\Spec(\scrR)$ is an affine, integral scheme such that the following three properties hold: 

\medskip
{\bf (i)} if $(\bar N,\phi_{\bar N},\vartheta_{\bar N},\nabla_{\bar N})$ is the evaluation of $\dbD(\scrD)$ at the trivial thickening $\grS_1(\scrR)$, then the $\scrR$-module $\bar N$ is free of rank $r$;

\smallskip
{\bf (ii)} the kernel $F^1_{\bar N}$ of the $\sigma_{\scrR}$-linear endomorphism $\phi_{\bar N}:\bar N\to\bar N$ is a direct summand of $\bar N$ which is a free $\scrR$-module of rank $d$;

\smallskip
{\bf (iii)} there exists an \'etale $k$-monomorphism $k[x_1,\ldots,x_{cd}]\hookrightarrow\scrR$ such that each $x_i$ is mapped to the maximal ideal $\scrL_0$ of $\scrR$ that defines the point $y_0$.

\medskip\noindent
Above $\vartheta_{\bar N}$ is the Verschiebung map of $\Phi_{\bar N}$ and $\nabla_{\bar N}$ is a connection on $\bar N$. As $\scrD$ is a versal deformation at all $k$-valued points of $\scrA$ and as $\scrA$ has dimension $cd$, we get:

\medskip
{\bf (iv)} the Kodaira--Spencer map $\grk_{\bar N}$ of $\nabla_{\bar N}$ is an $\scrR$-linear isomorphism. 

\medskip
Let $\scrR^{\text{l}}$ be a $p$-adically complete, formally smooth $W(k)[x_1,\ldots,x_{cd}]$-algebra which modulo $p$ is the $k[x_1,\ldots,x_{cd}]$-algebra $\scrR$. Let $\Phi_{\scrR^{\text{l}}}$ be the Frobenius lift of $\scrR^{\text{l}}$ that is compatible with $\sigma$ and that takes $x_i$ to $x_i^p$ for all $i\in\{1,\ldots,cd\}$. Let $\break d\Phi_{\scrR^{\text{l}}}:\Omega_{\scrR^{\text{l}}}^\wedge\to\Omega_{\scrR^{\text{l}}}^\wedge$ be the differential map of $\Phi_{\scrR^{\text{l}}}$. Let $(N,\phi_N,\vartheta_N,\nabla_N)$ be the projective limit indexed by positive integers $l$ of the evaluations of $\dbD(\scrD)$ at the thickenings attached naturally to the closed embeddings $\Spec(\scrR)\hookrightarrow\Spec(\scrR^{\text{l}}/p^l\scrR^{\text{l}})$; its reduction modulo $p$ is $(\bar N,\phi_{\bar N},\vartheta_{\bar N},\nabla_{\bar N})$. Due to the property (i), the $\scrR^{\text{l}}$-module $N$ is free of rank $r$. Let $F^1_N$ be a direct summand of $N$ that lifts $\bar F^1_N$. Based on properties (i) and (ii), there exists an isomorphism $(N,F^1_N)\arrowsim (M\otimes_{W(k)} \scrR^{\text{l}},F^1\otimes_{W(k)} \scrR^{\text{l}})$ to be viewed as a (non-canonical) identification. Under this identification,  $\phi_N$ and $\nabla_N$ get identified with $g_{\scrA}(\phi\otimes\Phi_{\scrR^{\text{l}}})$ for some element $g_{\scrA}\in \pmb{GL}_M(\scrR^{\text{l}})$ and with a connection $\nabla_M$ on $M\otimes_{W(k)} \scrR^{\text{l}}$ (respectively). We have: 
$$\nabla_M\circ g_{\scrA}(\phi\otimes\Phi_{\scrR^{\text{l}}})=(g_{\scrA}(\phi\otimes\Phi_{\scrR^{\text{l}}})\otimes d\Phi_{\scrR^{\text{l}}})\circ \nabla_M:(M+{1\over p}F^1)\otimes_{W(k)} \scrR^{\text{l}}\to M\otimes_{W(k)} \Omega_{\scrR^{\text{l}}}^\wedge.\leqno (7a)$$  
The Kodaira--Spencer map of the reduction modulo $p$ of $\nabla_M$ is an $\scrR$-linear map
$$\bar\grK_M:\oplus_{i=1}^{cd} \scrR{{\partial}\over {\partial x_i}}\to \Hom(\bar F^1,\bar F^0)\otimes_k \scrR=[\Lie(\pmb{GL}_{M/pM})/\Lie(\scrW_{+0k})]\otimes_k \scrR\leqno (7b)$$ 
which (due to the property (iv)) is an isomorphism. 

Let $\scrL_0^{\text{l}}$ be the ideal of $\scrR^{\text{l}}$ such that we have $\scrR^{\text{l}}/\scrL_0^{\text{l}}=\scrR/\scrL_0=k$. The reduction of ${1\over p}d\Phi_{\scrR^{\text{l}}}$ modulo $\scrL_0^{\text{l}}$ is $0$. This implies that the reduction modulo $\scrL_0^{\text{l}}$ of the following map $(g_{\scrA}(\phi\otimes\Phi_{\scrR^{\text{l}}})\otimes d\Phi_{\scrR^{\text{l}}})\circ \nabla_M:(M+{1\over p}F^1)\otimes_{W(k)} \scrR^{\text{l}}\to M\otimes_{W(k)} \Omega_{\scrR^{\text{l}}}^\wedge$ is $0$. Thus the reduction modulo $\scrL_0^{\text{l}}$ of the map $\nabla_M\circ g_{\scrA}(\phi\otimes\Phi_{\scrR^{\text{l}}}):(M+{1\over p}F^1)\otimes_{W(k)} \scrR^{\text{l}}\to M\otimes_{W(k)} \Omega_{\scrR^{\text{l}}}^\wedge$ is also $0$, cf. (7a). Therefore $\nabla_M$ modulo $\scrL_0^{\text{l}}$ is $\delta_0-g_{\scrA}^{-1}dg_{\scrA}$ modulo $\scrL_0^{\text{l}}$. Due to this and (7b),  we get that by replacing $\scrA=\Spec(\scrR)$ with an affine, open subscheme of it that contains the point $y_0\in\scrA(k)$, we can assume that: 

\medskip
{\bf (v)} the composite of the morphism $g_{\scrA}[1]:\scrA\to\pmb{GL}_{M/pM}$ defined by $g_{\scrA}$ modulo $p$ with the natural quotient morphism $\pmb{GL}_{M/pM}\twoheadrightarrow \pmb{GL}_{M/pM}/\scrW_{+0k}$, is an \'etale morphism.

\medskip
We consider the unique $W(k)$-monomorphism $\nu:\scrR^{\text{l}}\hookrightarrow W(\scrR)$ that lifts the identification $\scrR^{\text{l}}/p\scrR^{\text{l}}=\scrR$ and that takes $x_i\in\scrR^{\text{l}}$ to $(x_i,0,\ldots)\in W(\scrR)$ for all $i\in\{1,\ldots,cd\}$. The $W(k)$-homomorphism $\nu$ is compatible with the Frobenius lifts and it allows us to view $\pmb{GL}_M(\scrR^{\text{l}})$ as a subgroup of  $\pmb{GL}_M(W(\scrR))$. Thus we have $g_{\scrA}\in \pmb{GL}_M(W(\scrR))$. Let $g_{\scrA}[m]\in\pmb{GL}_M(W_m(\scrR))=\scrD_m(\scrR)$ be the natural reduction of $g_{\scrA}$ and let  
$$\eta_m:\scrA\to\scrD_m$$ 
be the morphism defined by $g_{\scrA}[m]$. 

\medskip\noindent
{\bf 3.1.1. End of the proof of 1.2 (a).} Let $y\in\scrA(k)$. Let $g_y\in \pmb{GL}_M(W(k))$ be the pull back of $g_{\scrA}$ via the Teichm\"uller section $\Spec(W(k))\hookrightarrow \Spec(W(\scrR))$ defined by $y$. The Dieudonn\'e module of $y^*(\scrD)$ is $(M,g_y\phi)$. The triple $(M/p^mM,g_y[m]\phi_m,\vartheta_m g_y[m]^{-1})$ is isomorphic to $(M/p^mM,\phi_m,\vartheta_m)$ (i.e., $y^*(\scrD)[p^m]$ is isomorphic to $D[p^m]$) if and only if we have $\eta_m(k)(y)\in\scrO_m(k)\subseteq\scrD_m(k)$, cf. Lemma 2.2.1. Therefore we have $y\in\grs_D(m)(k)$  if and only if $\eta_m(k)(y)\in\scrO_m(k)$. Thus we can take $\grs_D(m)$ to be the reduced, locally closed subscheme of $\scrA$ which is the reduced scheme of $\eta_m^*(\scrO_m)$.
\endproof

\medskip\noindent
{\bf 3.1.2. Remark.} Let $k_1$ be an algebraically closed field that contains $k$. The action $\dbT_m\times_k k_1$ over $k_1$ is the analogue action associated to the $p$-divisible group $D_{k_1}$ instead of to $D$. From this and the definition of $\grs_D(m)$ in Subsection 3.1.1, we get that we have an identity $\grs_D(m)(k_1)=\{y\in\scrA(k_1)|y^*(\scrD)[p^m]\;\text{is}\;\text{isomorphic}\;\text{to}\;D_{k_1}[p^m]\}$. In other words, the strata $\grs_D(m)$ are compatible with pulls back via geometric points. 

\medskip\smallskip\noindent
{\bf 3.2. Proof of 1.2 (b).} Let $y_1,y_2\in\grs_D(m)(k)$. For $j\in\{1,2\}$, let $I_j$ be the completion of the local ring of $\scrA$ at $y_j$. Let $\grs_j:=\Spec(I_j)\times_{\scrA} \grs_D(m)$; it is a reduced, local, complete, closed subscheme of $\Spec(I_j)$. Let $\grD_j$ be the pull back of $\scrD$ via the natural formally \'etale morphism $\Spec(I_j)\to\scrA$. As $\scrA$ has dimension $cd$ and as $\scrD$ is a versal deformation at all $k$-valued points of $\scrA$, $\grD_j[p^m]$ is the universal deformation of $y_j^*(\scrD)[p^m]$ (cf. [Il, Cor. 4.8 (ii)]). We can identify (non-canonically) $y_1^*(\scrD)[p^m]=D[p^m]=y_2^*(\scrD)[p^m]$. From the last two sentences, we get that there exists a unique isomorphism $\gamma_{12}:\Spec(I_1)\arrowsim\Spec(I_2)$ for which we have a unique isomorphism $\grD_1[p^m]\arrowsim \gamma_{12}^*(\grD_2)[p^m]$ that lifts the identification $y_1^*(\scrD)[p^m]=y_2^*(\scrD)[p^m]$. Based on Remark 3.1.2 we easily get that we have an identity
$$\grs_1=\gamma_{12}^*(\grs_2)$$
of reduced schemes. This identity implies that:

\medskip
{\bf (i)} the $k$-scheme $\grs_1$ is regular if and only if the $k$-scheme $\grs_2$ is regular, and

\smallskip
{\bf (ii)} we have $\dim(\grs_1)=\dim(\grs_2)$.

\medskip
Let $\grr_D(m)$ be a connected, open, regular subscheme of $\grs_D(m)$. We take $y_1$ such that it is a $k$-valued point of $\grr_D(m)$. From the property (i) we get that $\grs_D(m)$ is regular at any other $k$-valued point $y_2$ of it. Thus $\grs_D(m)$ is a regular scheme. From this and the property (ii) we get that all local rings of $\grs_D(m)$ of residue field $k$, have dimension $\dim(\grr_D(m))$. Thus $\grs_D(m)$ is also equidimensional.\endproof

\medskip\smallskip\noindent
{\bf 3.3. Proof of 1.2 (c).} To prove Theorem 1.2 (c), we can work locally in the Zariski topology of $\scrA$ and therefore we can assume that we are in the context described in Subsection 3.1. Let the morphism $\eta_m:\scrA\to\scrD_m$ be as in Subsection 3.1. We consider the composite epimorphism (induced by $\text{Red}_{l,\pmb{GL}_M}$ with $2\le l\le m$)
$$\pi_m:\scrD_m\twoheadrightarrow\scrD_{m-1}\twoheadrightarrow\cdots\twoheadrightarrow\scrD_1=\pmb{GL}_{M/pM}\twoheadrightarrow \pmb{GL}_{M/pM}/\scrW_{+0k}.$$ 
The composite morphism $\pi_m\circ\eta_m$ is \'etale (cf. property 3.1 (v)) and therefore its image $\scrU_m$ is an open subscheme of $\pmb{GL}_{M/pM}/\scrW_{+0k}$. The morphism $\eta_m$ is a section in the \'etale topology of $\scrU_m$ of the epimorphism $\scrD_m\times_{\pmb{GL}_{M/pM}/\scrW_{+0k}} \scrU_m\twoheadrightarrow\scrU_m$ induced naturally by $\pi_m$. This implies that $\scrB:=\im(\eta_m)$ and $\grw_D(m):=\im(\eta_m(\grs_D(m)))$ are reduced, locally closed subscheme of $\scrD_m$ and $\scrO_m$ (respectively). The natural morphisms $\scrA\to\scrB$ and $\grs_D(m)\to\grw_D(m)$ are \'etale. Thus $\scrB$ and $\grw_D(m)$ are regular $k$-schemes which are equidimensional of dimensions $cd$ and $\dim(\grs_D(m))$ (respectively).

For $\diamond$ a subscheme of $\scrD_m$, let $\dbT_{m,\diamond}:\scrH_m\times_k\diamond\to\scrD_m$ be the restriction of the action $\dbT_m$ to $\scrH_m\times_k\diamond$. Let $\scrE:=\dbT_{m,\scrB}^*(1_M[m])$; it is a closed subscheme of $\scrH_m\times_k\scrB$. In the next four paragraphs we check that the following property holds:

\medskip
{\bf (*)} the morphism $\dbT_{m,\scrB}:\scrH_m\times_k \scrB\to\scrD_m$ is smooth at all $k$-valued points of $\scrE$. 

\medskip
As $\dbT_{m,\scrB}:\scrH_m\times_k \scrB\to\scrD_m$ is a morphism between smooth $k$-schemes, to check the property (*) it suffices to show that for each point $(h_0[m],g_0[m])\in\scrE(k)$, the tangent map 
$$d\dbT_{m,\scrB}^{(h_0[m],g_0[m])}:L_{h_0[m]}\oplus L_{g_0[m]}\to\Lie(\scrD_m)$$ at the point $(h_0[m],g_0[m])\in\scrE(k)$ is onto. Here $L_{h_0[m]}$ and $L_{g_0[m]}$ are the tangent spaces of $\scrH_m$ and $\scrB$ (respectively) at its $k$-valued points $h_0[m]$ and $g_0[m]$ (respectively). Let $l_{h_0[m]}:\scrD_m\arrowsim\scrD_m$ be the isomorphism defined by the left translation via the action $\dbT_m$ through the point $h_0[m]\in\scrH_m(k)$. Let $r_{h_0[m]}:\scrH_m\arrowsim\scrH_m$ be the right translation through $h_0[m]\in\scrH_m(k)$. The isomorphism $r_{h_0[m]}\times l_{h_0[m]}^{-1}:\scrH_m\times_k l_{h_0[m]}(\scrB)\arrowsim \scrH_m\times_k \scrB$ is such that we have an identity $\dbT_{m,\scrB}\circ r_{h_0[m]}\times l_{h_0[m]}^{-1}=\dbT_{m,l_{h_0[m]}(\scrB)}:\scrH_m\times_k l_{h_0[m]}(\scrB)\to\scrD_m$. The  isomorphism $r_{h_0[m]}\times l_{h_0[m]}^{-1}$ maps  the point $((1_M[m],1_M[m],1_M[m]),1_M[m])\in \scrH_m(k)\times l_{h_0[m]}(\scrB)(k)$ to $(h_0[m],g_0[m])\in\scrH_m(k)\times \scrB(k)$; this is so as the relation $(h_0[m],g_0[m])\in\scrE(k)$ implies that $\dbT_m(h_0[m],g_0[m])=1_M[m]$. Thus by replacing the morphism $\eta_m$ with $l_{h_0[m]}\circ\eta_m$, to check that the tangent map $d\dbT_{m,\scrB}^{(h_0[m],g_0[m])}$ is onto, we can assume that $h_0[m]=(1_M[m],1_M[m],1_M[m])$ and $g_0[m]=1_M[m]$; thus we have $L_{h_0[m]}=\Lie(\scrH_m)$. 

We check by induction on $m\ge 1$ that $d\dbT_{m,\scrB}^{(h_0[m],1_M[m])}$ is onto. As $\dbT_m(h[m],1_M[m])$ is the product of the elements $h_1[m]h_2[m]h_3[m]^p$ and $\sigma_{\phi}(h_3[m])^{-1}\sigma_{\phi}(h_2[m])^{-1}\sigma_{\phi}(h_1[m]^p)^{-1}$ of $\scrD_m(k)=\pmb{GL}(W_m(k))$ (cf. (1a)) and as the Frobenius endomorphism of the $k$-algebra $k\oplus k\eps$ with $\eps^2=0$ annihilates the ideal $k\eps$, the restriction of $d\dbT_{m,\scrB}^{(h_0[m],1_M[m])}$ to $L_{h_0[m]}=\Lie(\scrH_m)$ is the same as the differential map of the homomorphism
$\dbW_m(\scrP_{0-}):\scrH_m\to\scrD_m$ which is defined by the homomorphism $\scrP_{0-W_m(k)}=\tilde\scrP_{0-W_m(k)}:\scrH_{W_m(k)}=\tilde\scrH_{W_m(k)}\to\pmb{GL}_{M/p^mM}$ (see Subsubsection 2.1.1) and which maps the point $(h_1[m],h_2[m],h_3[m])\in\scrH_m(k)$ to $h_1[m]h_2[m]h_3[m]^p\in\scrD_m(k)=\pmb{GL}(W_m(k))$. 

Suppose that $m=1$. We identify $\Lie(\scrH_1)=\Lie(\scrW_{+k})\oplus\Lie(\scrW_{0k})\oplus\Lie(\scrW_{-k})$ as $k$-vector spaces. The restriction of $d\dbT_{1,\scrB}^{(h_0[1],1_M[1])}$ to the direct summand $\Lie(\scrW_{+k})\oplus\Lie(\scrW_{0k})$ of $\Lie(\scrH_1)\oplus L_{1_M[1]}$ is injective and its image is the Lie subalgebra $\Lie(\scrW_{+0k})$ of $\Lie(\scrD_1)=\Lie(\pmb{GL}_{M/pM})$, cf. the definition of $\dbW_1(\scrP_{0-})$. The $k$-vector space $d\dbT_{1,\scrB}^{(h_0[1],1_M[1])}(L_{1_M[1]})$ is a direct supplement of $\Lie(\scrW_{+0k})$ in $\Lie(\scrD_1)=\Lie(\pmb{GL}_{M/pM})$, cf. property 3.1 (v). From the last two sentences we get that $d\dbT_{1,\scrB}^{(h_0[1],1_M[1])}$ is onto. Thus the basis of the induction holds. For $m\ge 2$, the passage from $m-1$ to $m$ goes as follows. 

The action $\dbT_m$ is a natural lift of the action $\dbT_{m-1}$ and we have reduction epimorphisms $\text{Red}_{m,\tilde\scrH}:\scrH_m\twoheadrightarrow\scrH_{m-1}$ and $\text{Red}_{m,\pmb{GL}_M}:\scrD_m\twoheadrightarrow\scrD_{m-1}$ whose kernels are the vector groups over $k$ defined by $\Lie(\scrH_k)$ and $\Lie(\pmb{GL}_{M/pM})$ (respectively). Thus we have natural short exact sequences of Lie algebras 
$$0\to \Lie(\scrH_k)_{\text{ab}}\to\Lie(\scrH_m)\to\Lie(\scrH_{m-1})\to 0\leqno (8a)$$
and
$$0\to \Lie(\pmb{GL}_{M/pM})_{\text{ab}}\to\Lie(\scrD_m)\to\Lie(\scrD_{m-1})\to 0,\leqno (8b)$$
where  $\bot_{\text{ab}}$ is the abelian Lie algebra on the $k$-vector space $\bot$.
As $d\dbT_{m-1,\scrB}^{(h_0[m-1],1_M[m-1])}$ is onto, to check that $d\dbT_{m,\scrB}^{(h_0[m],1_M[m])}$ is onto it suffices to show that we have an inclusion $\Lie(\pmb{GL}_{M/pM})_{\text{ab}}\subseteq\im(d\dbT_{m,\scrB}^{(h_0[m],1_M[m])})$. But $d\dbT_{m,\scrB}^{(h_0[m],1_M[m])}(\Lie(\scrH_k)_{\text{ab}})$ is the Lie subalgebra $\Lie(\scrW_{+0k})_{\text{ab}}$ of $\Lie(\pmb{GL}_{M/pM})_{\text{ab}}$, cf. the definition of $\dbW_m(\scrP_{0-})$. Using a smooth parametric curve which has the form $(1_M[m],1_M[m],h_3(t)[m])\in\scrH_m(k)$ for $t\in k$, where the element $h_3(t)[m]\in\scrW_-(W_m(k))$ is such that $h_3(0)[m]=1_M[m]$ and $h_3(t)[m]^p\in\Ker(\scrW_-(W_m(k))\to\scrW_-(W_{m-1}(k)))$, from the definition of $\dbW_m(\scrP_{0-})$ we get that $d\dbT_{m,\scrB}^{(h_0[m],1_M[m])}(\Lie(\scrH_m))$ contains the Lie subalgebra $\Lie(\scrW_{-k})_{\text{ab}}$ of $\Lie(\pmb{GL}_{M/pM})_{\text{ab}}$. As $\Lie(\pmb{GL}_{M/pM})_{\text{ab}}=\Lie(\scrW_{+0k})_{\text{ab}}\oplus\Lie(\scrW_{-k})_{\text{ab}}$, we get that the desired inclusion $\break\Lie(\pmb{GL}_{M/pM})_{\text{ab}}\subseteq\im(d\dbT_{m,\scrB}^{(h_0[m],1_M[m])})$ holds. Therefore $d\dbT_{m,\scrB}^{(h_0[m],1_M[m])}$ is onto. This ends the induction and thus also the proof of the property (*).

Due to the property (*), the $k$-scheme $\scrE$ is smooth of dimension equal to $\break\dim(\scrH_m)+\dim(\scrB)-\dim(\scrD_m)=\dim(\scrB)$. The reduced subscheme of the image of the natural projection morphism $\Pi_2:\scrE\to\scrB$ is $\grw_D(m)$ and the fibres of $\Pi_2$ above $k$-valued points of $\grw_D(m)$ are naturally identified with stabilizer subgroup schemes of $k$-valued points of $\grw_D(m)\subseteq \scrO_m$. Thus all the fibres of $\Pi_2$ above $k$-valued points of $\grw_D(m)$ have dimension $\gamma_D(m)$, cf. Theorem  2.4 (c). This implies that $\dim(\grw_D(m))=\dim(\scrE)-\gamma_D(m)$. Thus $\dim(\grw_D(m))=\dim(\scrB)-\gamma_D(m)$. We recall that $\dim(\scrB)=\dim(\scrA)=cd$ and that $\dim(\grw_D(m))=\dim(\grs_D(m))$. From the last two sentences we get the desired identity $\dim(\grs_D(m))=cd-\gamma_D(m)$.
\endproof

\medskip\noindent
{\bf 3.3.1. Corollary.} {\it For all $m\ge 1$ we have 
$$\dim(\scrO_m)=mr^2-cd+\dim(\grs_D(m)).\leqno (9)$$}
\noindent
{\it Proof:} We have $\dim(\scrO_m)=mr^2-\gamma_D(m)$, cf. Formula (1) and Theorem 2.4 (c). From this and Theorem 1.2 (c) we get that Formula (9) holds.\endproof 

\medskip\smallskip\noindent
{\bf 3.4. Proof of 1.2 (d).} If $y\in\grs_D(n_D)(k)$ (i.e., if $y^*(\scrD)[p^{n_D}]$ is isomorphic to $D[p^{n_D}]$), then $y^*(\scrD)$ is isomorphic to $D$ (cf. the very definition of $n_D$). Thus for all $m\ge n_D$, $y^*(\scrD)[p^m]$ is isomorphic to $D[p^m]$ i.e., we also have $y\in\grs_D(m)(k)$. Therefore for all $m\ge n_D$, we have $\grs_D(n_D)\subseteq\grs_D(m)$; as we obviously have $\grs_D(m)\subseteq\grs_D(n_D)$, we conclude that $\grs_D(m)=\grs_D(n_D)$. Based on this and Theorem 1.2 (c) we get that for all $m\ge n_D$ we have $\gamma_D(m)=s_D$. 
\endproof

\medskip\smallskip\noindent
{\bf 3.5. Proof of 1.2 (e).} Let $D\to \tilde D$ be an isogeny of $p$-divisible group over $k$. Let $\kappa\in\dbN$ be such that $p^{\kappa}$ annihilates the kernel $\scrK$ of this isogeny. We identify naturally $\tilde D=D/\scrK$. We take $m$ such that we have $m\ge\max\{\kappa+n_{\tilde D},n_D\}$. To prove that $s_D=s_{\tilde D}$, we can assume that there exists a point $y_{\tilde D}:\Spec(k)\to\scrA$ such that $y^*_{\tilde D}(\scrD)$ is isomorphic to $\tilde D$ (if needed, one can add extra connected components to  $\scrA$). Let $I_0$ be the completion of the local ring of $\scrA$ at $y_{\tilde D}$. Let $\grs_{\tilde D}(m-\kappa)$ be the reduced, locally closed subscheme of $\scrA$ such that we have an identity $\grs_{\tilde D}(m-\kappa)(k)=\{\tilde y\in\scrA(k)|\tilde y^*(\scrD)[p^{m-\kappa}]\;\text{is}\;\text{isomorphic}\;\text{to}\;\tilde D[p^{m-\kappa}]\}$ (cf. Theorem 1.2 (a)). We have $\dim(\grs_{\tilde D}(m-\kappa))=cd-s_{\tilde D}$ and $\dim(\grs_D(m))=cd-s_{D}$, cf. Theorems 1.2 (c) and (d) and the inequality $m\ge\max\{\kappa+n_{\tilde D},n_D\}$.

Let $\gri_D(m)$ be a finite, flat $\grs_D(m)$-scheme which is smooth over $k$ and for which there exists an isomorphism $\scrD[p^m]\times_{\scrA} \gri_D(m)\arrowsim D[p^m]\times_k \gri_D(m)$, to be viewed as a natural identification (cf. [Va2, Thm. 5.3.1 (c)]). Let $\tilde\scrD_{m,m}$ be the quotient of $\scrD\times_{\scrA} \gri_D(m)$ by its finite, flat subgroup scheme $\scrK\times_k \gri_D(m)$; it is a $p$-divisible group scheme over $\gri_D(m)$. For each closed point $\tilde y:\Spec(k)\hookrightarrow\gri_D(m)$, the pull back $\tilde y^*(\tilde\scrD_{m,m})$ is a $p$-divisible group over $k$ whose Barsotti--Tate group of level $m-\kappa$ is isomorphic to $\tilde D[p^{m-\kappa}]$. As $m-\kappa\ge n_{\tilde D}$, we get that $\tilde y^*(\tilde\scrD_{m,m})$ is isomorphic to $\tilde D$. If $\tilde I$ is the completion of the local ring of $\gri_D(m)$ at $\tilde y$, then the $p$-divisible group $\tilde\scrD_{m,m}\times_{\gri_D(m)} \Spec(\tilde I)$ over $\Spec(\tilde I)$ is the pull back of $\scrD\times_{\scrA} \Spec(I_0)$ via a composite morphism $\Spec(\tilde I)\to \grs_{\tilde D}(m-\kappa)\times_{\scrA} \Spec(I_0)\hookrightarrow \Spec(I_0)$. We show that the assumption that the morphism $\Spec(\tilde I)\to \grs_{\tilde D}(m-\kappa)\times_{\scrA} \Spec(I_0)$ does not have a finite fibre over $k$ leads to a contradiction. This assumption implies that there exists an integral, closed subscheme $\tilde Y$ of $\Spec(\tilde I)$ which is of positive dimension and over which the natural pull back of $\scrD$ is constant (i.e., it is isomorphic to $\tilde D\times_k \tilde Y$). This implies that the natural morphism $\tilde Y\to\scrA$ is constant (i.e., it factors through $\Spec(k)$). As the morphism $\gri_D(m)\to\scrA$ is quasi-finite, we get that $\tilde Y$ has dimension $0$. Contradiction. 

As the morphism $\Spec(\tilde I)\to \grs_{\tilde D}(m-\kappa)\times_{\scrA} \Spec(I_0)$ has a finite fibre over $k$, at the level of dimensions we have the following relations
$$cd-s_D=\dim(\grs_D(m))=\dim(\gri_D(m))=\dim(\Spec(\tilde I))$$
$$\le \dim(\grs_{\tilde D}(m-\kappa)\times_{\scrA} \Spec(I_0))=\dim(\grs_{\tilde D}(m-\kappa))=cd-s_{\tilde D}.$$
Thus $s_{\tilde D}\le s_D$. Interchanging the roles of $D$ and $\tilde D$, a similar argument shows that $s_D\le s_{\tilde D}$. Thus $s_D=s_{\tilde D}$ i.e., the specializing height $s_D$ of $D$ is an isogeny invariant. 
\endproof

\medskip\smallskip\noindent
{\bf 3.6. Proof of 1.2 (f).} If $E_1$ and $E_2$ are two finite, commutative group schemes over $k$, then $\pmb{\text{Hom}}(E_1,E_2)$ is the affine group scheme of finite type over $k$ that parametrizes homomorphisms between $E_1$ and $E_2$. The scheme $\pmb{\text{Aut}}(E_1)$ is a non-empty, open subscheme of $\pmb{\text{End}}(E_1)$. As $\pmb{\text{End}}(E_1)$ is equidimensional, we conclude that:

\medskip
{\bf (i)} we have $\dim(\pmb{\text{Aut}}(E_1))=\dim(\pmb{\text{End}}(E_1))$. 
\medskip

Based on Theorem 1.2 (c), to prove Theorem 1.2 (f) it suffices to prove the following identity $s_D=cd-{1\over 2}\sum_{s=1}^v\sum_{t=1}^v r_sr_t|\alpha_s-\alpha_t|$. As both numbers $cd-{1\over 2}\sum_{s=1}^v\sum_{t=1}^v r_sr_t|\alpha_s-\alpha_t|$ and $s_D$ are isogeny invariants (cf. Theorem 1.2 (e) for $s_D$), to check that they are equal we can replace $D$ by any other $p$-divisible group over $k$ isogenous to it (this idea appears first in [Tr2, Sect. 1, pp. 46--47]). Thus we can assume that $D$ is {\it minimal} in the sense of [Oo3, Subsect. 1.1] (this idea appears first in an informal manuscript of Oort and it slightly shortens the computations; we recall that Traverso used in [Tr2, Sect. 1, pp. 46--47] direct sums of special $p$-divisible groups over $k$ of $a$-number at most $1$). Thus we have a product decomposition $D=\prod_{s=1}^v D_s$, where $D_s$ is a minimal $p$-divisible group of height $r_s$ and Newton polygon slope $\alpha_s={{d_s}\over {r_s}}$. To simplify the calculations, we can assume that $\alpha_1\le\alpha_2\le\cdots\le\alpha_m$. If $1\le t\le s\le v$, then $r_sr_t|\alpha_s-\alpha_t|=c_sd_t-c_td_s$. Thus we compute
$$cd-{1\over 2}\sum_{s=1}^v\sum_{t=1}^v r_sr_t|\alpha_s-\alpha_t|=cd-\sum_{s=2}^v\sum_{t=1}^{s-1} (c_sd_t-c_td_s)\leqno (10a)$$
$$=\sum_{s=1}^v c_sd_s+\sum_{s=2}^v\sum_{t=1}^{s-1} (c_sd_t+c_td_s)-\sum_{s=2}^v\sum_{t=1}^{s-1} (c_sd_t-c_td_s)=\sum_{s=1}^v c_sd_s+2\sum_{s=2}^v\sum_{t=1}^{s-1} c_td_s.$$ 
\indent
We have $\dim(\pmb{\text{Aut}}(D_s[p]))=c_sd_s$, cf. Example 2.3.1 and Theorem 2.4 (c) applied to $(D_s,c_s,d_s)$ instead of to $(D,c,d)$. Thus $\dim(\pmb{\text{End}}(D_s[p]))=c_sd_s$, cf. property (i).

If $1\le t<s\le m$ and $\alpha_s<\alpha_t$, then using the property (i) we have
$$\dim(\pmb{\text{Aut}}(D_s[p]\times_k D_t[p]))=\dim(\pmb{\text{End}}(D_s[p]\times_k D_t[p]))=\dim(\pmb{\text{End}}(D_s[p]))+\dim(\pmb{\text{End}}(D_t[p]))$$
$$+\dim(\pmb{\text{Hom}}(D_t[p],D_s[p]))+\dim(\pmb{\text{Hom}}(D_s[p],D_t[p]))=c_sd_s+c_td_t+2c_td_s$$ 
(cf. end of Example 2.3.2 and Theorem 2.4 (c) applied to $(D_s\times_k D_t,c_s,c_t,d_s,d_t)$). Thus $\dim(\pmb{\text{Hom}}(D_t[p],D_s[p]))+\dim(\pmb{\text{Hom}}(D_s[p],D_t[p]))=2c_td_s$ (either the third paragraph of Example 2.3.2 or standard Cartier duality can be used to show that in fact we have $\dim(\pmb{\text{Hom}}(D_t[p],D_s[p]))=\dim(\pmb{\text{Hom}}(D_s[p],D_t[p]))=c_td_s$). The identities of this paragraph also hold if $\alpha_s=\alpha_t$ (the references to Example 2.3.2 being replaced by references to Example 2.3.1).

We have $n_D=1$, cf. either [Oo3, Thm. 1.2] or [Va4, Main Thm. D]. Thus 
$$s_D=\gamma_D(1)=\dim(\pmb{\text{Aut}}(D[p]))=\dim(\pmb{\text{End}}(D[p]))=\sum_{s=1}^v\dim(\pmb{\text{End}}(D_s[p]))\leqno (10b)$$
$$+\sum_{s=2}^v\sum_{t=1}^{s-1} [\dim(\pmb{\text{Hom}}(D_t[p],D_s[p]))+\dim(\pmb{\text{Hom}}(D_s[p],D_t[p]))]=\sum_{s=1}^v c_sd_s+2\sum_{s=2}^v\sum_{t=1}^{s-1} c_td_s.$$
From Formulas (10a) and (10b) we get that $s_D=cd-{1\over 2}\sum_{s=1}^v\sum_{t=1}^v r_sr_t|\alpha_s-\alpha_t|$.\endproof

\medskip\smallskip\noindent
{\bf 3.7. Proofs of 1.2 (g) and (h).} To prove Theorems 1.2 (g) and (h), we can work locally in the Zariski topology of $\scrA$ and therefore we can assume that we are in the context described in Subsection 3.1; in particular, $\scrA=\Spec(\scrR)$ is affine. As $\scrO_m$ is a quasi-affine scheme (see paragraph before Lemma 2.2.10), $\eta_m^*(\scrO_m)$ is a quasi-affine $\scrA$-scheme. From this and Subsubsection 3.1.1, we get that $\grs_D(m)=(\eta_m^*(\scrO_m))_{\text{red}}$ is a quasi-affine $\scrA$-scheme. 

In order to prove a stronger form of Theorem 1.2 (h), in this paragraph we will only assume that the image of the homomorphism $\chi_D(m):\pmb{\text{Aut}}(D[p^m])(k)\to\scrW_0(k)$ introduced in Lemma 2.5, has finite image. This assumption implies that the orbit $\scrO_m$ is an affine scheme, cf. Lemmas 2.5 and 2.6. Thus the scheme $\grs_D(m)=(\eta_m^*(\scrO_m))_{\text{red}}$ is affine and therefore it is also an affine $\scrA$-scheme. This ends the proofs of Theorems 1.2 (g) and (h). This ends the proof of the Basic Theorem.\endproof

\medskip\noindent
{\bf 3.7.1. Remarks.} {\bf (a)} The computations of Subsection 2.3 that led to the proof of properties 2.3 (viii) and (ix), can be adapted to show that the reduced group of $\Ker(\scrC_m\to\scrC_{m-1})$ does not depend on $m\ge 2$ and moreover has dimension $\gamma_D(1)$. Thus if we have $\gamma_D(2)>\gamma_D(1)$ (for instance, this holds if $n_D=2$), then the group $\im(\scrC_2\to\scrC_1)$ (equivalently, the group $\im(\pmb{\text{Aut}}(D[p^2])\to\pmb{\text{Aut}}(D[p]))$ cf. Theorem 2.4 (b)) has positive dimension $\gamma_D(2)-\gamma_D(1)$ and therefore Theorem 1.2 (h) does not apply for $m=2$. Moreover, easy examples show that the image of the homomorphism $\chi_D(2)$ can be infinite and therefore in such cases even the previous paragraph does not apply for $m=2$. 

\smallskip
{\bf (b)} In practice, it is easy to decide if the homomorphism $\chi_D(m)$ does or does not have a finite image. But even if the homomorphism $\chi_D(m)$ has an infinite image, one expects that in many cases $\grs_D(m)$ is an affine $\scrA$-scheme.

\bigskip
\noindent
{\boldsectionfont 4. Applications to Shimura varieties of Hodge type}
\bigskip
In this Section we  show how the Basic Theorem transfers naturally to the context of special fibres of good integral models in unramified mixed characteristic $(0,p)$ of Shimura varieties of Hodge type. To be short (i.e., in order not to recall all the machinery of Shimura varieties of Hodge type) and for the sake of generality, we will work in the (more general and) axiomatized context of quasi Shimura $p$-varieties of Hodge type. In Subsection 4.1 we present the relative version of the orbit spaces of Section 2. In Subsection 4.2 we introduce the mentioned axiomatized context. An analogue of the Basic Theorem for  quasi Shimura $p$-varieties of Hodge type is presented in the Basic Corollary  4.3. Example 4.4 pertains to the applicability of Corollary 4.3 (e). Example 4.5 pertains to special fibres of Mumford's moduli schemes $\scrA_{d,1,l}$. Example 4.6 pertains to good integral models in mixed characteristic $(0,p)$ of Shimura varieties of Hodge type. We will use the notations listed before Subsection 2.1; thus $m$ is a positive integer. Once a good direct sum decomposition $M=F^1\oplus F^0$ is introduced, we will also use the notations listed before Subsection 3.1. 

\medskip\smallskip\noindent
{\bf 4.1. Relative orbit spaces.} Let $G$ be a smooth, closed subgroup scheme of $\pmb{GL}_M$ such that its generic fibre $G_{B(k)}$ is connected. Thus the scheme $G$ is integral. Until the end we will assume that the following two axioms hold for the triple $(M,\phi,G)$:

\medskip
{\bf (i)} the Lie subalgebra $\Lie(G_{B(k)})$ of $G_{B(k)}$ is stable under $\phi$ i.e., we have $\break\phi(\Lie(G_{B(k)}))=\Lie(G_{B(k)})$;

\smallskip
{\bf (ii)} there exist a direct sum decomposition $M=F^1\oplus F^0$ and a smooth, closed subgroup scheme $G_1$ of $\pmb{GL}_M$ such that the following four properties hold:

\medskip
\item{\bf (ii.a)} the kernel of the reduction modulo $p$ of $\phi$ is $F^1/pF^1$;

\smallskip
\item{\bf (ii.b)} the cocharacter $\mu:\dbG_m\to\pmb{GL}_M$ which acts trivially on $F^0$ and via the inverse of the identical character of $\dbG_m$ on $F^1$, factors through $G_1$;

\smallskip
\item{\bf (ii.c)} the group scheme $G$ is a normal, closed subgroup scheme of $G_1$ and we have a short exact sequence $1\to G\to G_1\to\dbG_m^u\to 1$, where $u\in\{0,1\}$;

\smallskip
\item{\bf (ii.d)} if $u=1$ (i.e., if $G_1\neq G$), then the homomorphism $\mu:\dbG_m\to G_1$ defined by $\mu$ (cf. property (ii.b)) is a splitting of the short exact sequence of the property (ii.c).

\medskip
If $G$ is a reductive group scheme and $u=0$, then the triple $(M,\phi,G)$ is called a {\it Shimura $F$-crystal} over $k$ (cf. [Va3, Subsect. 1.1]). In general, the triple $(M,\phi,G)$ is called an {\it $F$-crystal with a group} over $k$ (cf. [Va2, Def. 1.1 (a) and Subsect. 2.1]); if $u=0$, then the {\it $W$-condition} of [Va2, Subsubsect. 2.2.1 (d)] holds for $(M,\phi,G)$. Due to the properties (ii.b) and (ii.c) we have a direct sum decomposition
$$\Lie(G)=\oplus_{i=-1}^1\tilde F^i(\Lie(G))\leqno (11)$$
such that $\mu$ acts via inner conjugation on $\tilde F^i(\Lie(G))$ as the $-i$-th power of the identical character of $\dbG_m$. Let $e_+$, $e_0$, and $e_-$ be the ranks of $\tilde F^1(\Lie(G))$, $\tilde F^0(\Lie(G))$, and $\tilde F^{-1}(\Lie(G))$ (respectively). Let $d_G:=\dim(G_k)=\dim(G_{B(k)})$. Due to (11) we have
$$d_G=e_++e_0+e_-.$$
\noindent
{\bf 4.1.1. Relative group schemes.} We will use the notations listed before Subsection 3.1 for the direct sum decomposition $M=F^1\oplus F^0$ of the axiom 4.1 (ii). We consider the following five closed subgroup schemes
$\scrW_+^G:=\scrW_+\cap G$, $\scrW_0^G:=\scrW_0\cap G$, $\scrW_-^G:=\scrW_-\cap G$, $\scrW_{+0}^G:=\scrW_{+0}\cap G$, and $\scrW_{0-}^G:=\scrW_{0-}\cap G$ of $G$. Let 
$$\scrH^G:=\scrW_+^G\times_{W(k)} \scrW_0^G\times_{W(k)} \scrW_-^G;$$ 
it is a closed subscheme of $\scrH$ such that $\scrH^G_{W_m(k)}$ is a closed subgroup subscheme of $\scrH_{W_m(k)}=\tilde\scrH_{W_m(k)}$ (we recall that we view the isomorphism $\scrP_{W_m(k)}:\scrH_{W_m(k)}\arrowsim\tilde\scrH_{W_m(k)}$ of $\Spec(W_m(k))$-schemes as a natural identification).

The group schemes $\scrW_+^G$ and $\scrW_-^G$ over $\Spec(W(k))$ are isomorphic to $\dbG_a^{e_+}$ and $\dbG_a^{e_-}$ (respectively). More precisely, if $R$ is a commutative $W(k)$-algebra, then we have 
$$\scrW_+^G(R)=1_{M\otimes_{W(k)} R}+\tilde F^1(\Lie(G))\otimes_{W(k)} R\;\;\text{and}\;\;\scrW_-^G(R)=1_{M\otimes_{W(k)} R}+\tilde F^{-1}(\Lie(G))\otimes_{W(k)} R.$$ 
\indent
Let $\scrW_0^{G_1}:=\scrW_0\cap G_1$. The group scheme $\scrW_0^{G_1}$ is the centralizer of the torus $\im(\mu)$ in $G_1$ and therefore it is a smooth group scheme over $\Spec(W(k))$, cf. [DG, Vol. III, Exp. XIX, 2.2]. If $u=0$, then $\scrW_0^G=\scrW_0^{G_1}$. If $u=1$, then due to the property 4.1 (ii.d) we have a short exact sequence $1\to \scrW_0^G\to\scrW_0^{G_1}\to \dbG_m\to 1$ which splits; thus the group scheme $\scrW_0^{G_1}$ is isomorphic to the semidirect product of $\scrW_0^G$ and $\dbG_m$. We conclude that regardless of what $u$ is, the group scheme $\scrW_0^G$ is smooth over $\Spec(W(k))$.

The Lie algebras of $W_+^G$, $W_0^G$, and $W_-^G$ are $\tilde F^1(\Lie(G))$, $\tilde F^0(\Lie(G))$, and $\tilde F^{-1}(\Lie(G))$ (respectively). This implies that the relative dimension of $\scrW_0^G$ is $e_0$. The smooth, affine scheme $\scrH^G$ has relative dimension $d_G$ over $\Spec(W(k))$. The natural product morphism $\scrP_0^G:\scrH^G\to G$ is induced naturally by the open embedding $\scrP_0:\scrH\hookrightarrow\pmb{GL}_M$ and therefore it is also an open embedding.  Let $\scrP_-^G:=1_{\scrW_+^G}\times 1_{\scrW_0^G}\times p1_{\scrW_-^G}:\scrH^G\to\scrH^G$. The composite morphism $\scrP_{0-}^G:=\scrP_0^G\circ\scrP_-^G:\scrH^G\to G$ has the property that its reduction $\scrP_{0-W_m(k)}^G:\scrH^G_{W_m(k)}\to G_{W_m(k)}$ modulo $p^m$ is a homomorphism of affine group schemes over $\Spec(W_m(k))$ which is a restriction of the homomorphism $\scrP_{0-W_m(k)}:\scrH_{W_m(k)}\to\pmb{GL}_{M/p^mM}$ (see Subsection 2.1.1 for $\scrP_{0-}$). 

\medskip\noindent
{\bf 4.1.2. Lemma.} {\it The group scheme $\scrW_{+0}^G$  over $\Spec(W(k))$ is the semidirect product of $\scrW_+^G$ and $\scrW_0^G$. Thus $\scrW_{+0}^G$ is smooth.}

\medskip
\proof
We have an identity $\Lie(\scrW_{+0B(k)}^G)=\tilde F^1(\Lie(G))[{1\over p}]\oplus \tilde F^0(\Lie(G))[{1\over p}]$ as well as an inclusion $\Lie(\scrW_{+0k}^G)\subseteq\tilde F^1(\Lie(G))\otimes_{W(k)} k\oplus \tilde F^0(\Lie(G))\otimes_{W(k)} k$. Thus
$$\dim(\scrW_{+0B(k)}^G)=\dim_{B(k)}(\Lie(\scrW_{+0B(k)}^G))=e_++e_0\ge\dim_k(\Lie(\scrW_{+0k}^G))\ge\dim(\scrW_{+0k}^G).$$ 
The Zariski closure of $\scrW_{+0B(k)}^G$ in $G$ is contained in $\scrW_{+0}^G$ and therefore we have $\dim(\scrW_{+0k}^G)\ge\dim(\scrW_{+0B(k)}^G)$. From the last two sentences we get that $\dim_k(\Lie(\scrW_{+0k}^G))=\dim(\scrW_{+0k}^G)$. This implies that $\scrW_{+0k}^G$ is a smooth group over $k$. Thus the affine group scheme $\scrW_{+0}^G$  is smooth over $\Spec(W(k))$ if and only if the natural reduction homomorphism $\scrW_{+0}^G(W(k))\to\scrW_{+0}^G(k)$ is onto. 

The semidirect product $\tilde\scrW_{+0}^G$ of $\scrW_+^G$ and $\scrW_0^G$ is a smooth subgroup scheme of  $\scrW_{+0}^G$. To check that the group scheme $\scrW_{+0}^G$ is smooth over $\Spec(W(k))$ and equal to $\tilde\scrW_{+0}^G$, it suffices to show that the group $\scrW_{+0}^G(k)$ is equal to its subgroup $\tilde\scrW_{+0}^G(k)=\scrW_+^G(k)\scrW_0^G(k)$. 

Let $h_{12}[1]\in\scrW_{+0}^G(k)\leqslant \scrW_{+0}(k)$. As $\scrW_{+0}$ is the semidirect product of $\scrW_+$ and $\scrW_0$, there exist unique elements $h_1[1]\in\scrW_+(k)$ and $h_2[1]\in\scrW_0(k)$ such that $h_{12}[1]=h_1[1]h_2[1]$. As the cocharacter $\mu_k:\dbG_m\to G_{1k}$ normalizes $\scrW_{+0k}^G$, by considering conjugates of $h_{12}[1]$ through $k$-valued points of $\im(\mu_k)$, we easily get that we have $h_1[1],h_2[1]\in\scrW_{+0k}^G(k)$. This implies that $h_1[1]\in \scrW_+(k)\cap G(k)=\scrW_+^G(k)$ and $h_2[1]\in\scrW_0(k)\cap G(k)=\scrW_0^G(k)$. Therefore $h_{12}[1]=h_1[1]h_2[1]\in\tilde\scrW_{+0}^G(k)$ and thus $\scrW_{+0}^G(k)=\tilde\scrW_{+0}^G(k)$. Therefore $\scrW_{+0}^G=\tilde\scrW_{+0}^G$.\endproof

\medskip\noindent
{\bf 4.1.3. The relative action $\dbT_m^G$.} Let $\tilde\scrH^G$ be the dilatation of $G$ centered on the smooth subgroup $\scrW^G_{+0k}$ of $G_k$; it is a smooth, closed subgroup scheme of $\tilde\scrH$. As in Subsubsection 2.1.1 we argue that we have a natural morphism $\scrP^G:\scrH^G\to\tilde\scrH^G$ of $\Spec(W(k))$-schemes which gives birth to an isomorphism $\scrP_{W_m(k)}^G:\scrH^G_{W_m(k)}\arrowsim\tilde\scrH^G_{W_m(k)}$ of $\Spec(W_m(k))$-schemes, to be viewed as a natural identification. Obviously, the group schemes structures on $\scrH_{W_m(k)}^G$ induced via the identification $\scrP^G_{W_m(k)}$ or via the identification of $\scrH_{W_m(k)}^G$ with a closed subgroup scheme of $\scrH_{W_m(k)}$, are equal. Let $\scrH_m^G:=\dbW_m(\scrH^G)=\dbW_m(\tilde\scrH^G)$; it is a smooth, affine group of dimension $md_G$ which is connected if and only if $\scrH_k^G$ (equivalently $\scrW_{0k}^G$) is connected (cf. Subsubsection 2.1.4). Let $\scrD_m^G:=\dbW_m(G)$; it is a smooth, affine $k$-scheme of dimension $md_G$ which is connected if and only if $G_k$ is connected (cf. Subsubsection 2.1.4). 

As $\phi$ is a $\sigma$-linear automorphism of $M[{1\over p}]$, the group $\{\phi h\phi^{-1}|h\in G(B(k))\}$ is the group of $B(k)$-valued points of the unique connected subgroup of $\pmb{GL}_{M[{1\over p}]}$ that has $\phi(\Lie(G_{B(k)}))$ as its Lie algebra (see [Bo, Ch. II, Subsect. 7.1] for the uniqueness part). As $\phi(\Lie(G_{B(k)}))=\Lie(G_{B(k)})$ (cf. the axiom 4.1 (i)), we conclude that $\break\{\phi h\phi^{-1}|h\in G(B(k))\}=G(B(k))$. A similar argument shows that $\{\sigma_{\phi}h\sigma_{\phi}^{-1}|h\in G(W(k))\}=G(W(k))$ (see Subsection 2.2 for $\sigma_{\phi}$). For each $h\in \im(\scrP_{0-}^G(W(k)))$ we have $\phi h\phi^{-1}\in G(W(k))=G(B(k))\cap\pmb{GL}_M(W(k))$. Therefore we have a unique action
$$\dbT_m^G:\scrH_m^G\times_k\scrD_m^G\to\scrD_m^G$$
which is the natural restriction of the action $\dbT_m$ of Subsection 2.2. Let $\scrO_m^G$ be the orbit of $1_M[m]\in\scrD_m^G(k)$ under the action $\dbT_m^G$. Let $\scrS_m^G$ be the stabilizer subgroup scheme of the point $1_M[m]\in\scrD_m^G(k)$ under the action $\dbT_m^G$; we have $\scrS_m^G:=\scrS_m\cap\scrH_m^G$. Let $\scrC_m^G$ be the reduced group of $\scrS_m^G$ and let $\scrC_m^{0G}$ be the identity component of $\scrC_m^G$.

\medskip\noindent
{\bf 4.1.4. Definitions.} {\bf (a)} Let $g_1,g_2\in G(W(k))$. By an {\it inner isomorphism} between the two quadruples $(M/p^mM,g_1[m]\phi_m,\vartheta_m g_1[m]^{-1},G_{W_m(k)})$ and  $(M/p^mM,g_2[m]\phi_m,\vartheta_m g_2[m]^{-1},G_{W_m(k)})$, we mean an element $g_3[m]\in G(W_m(k))$ which is an isomorphism between the two Dieudonn\'e modules $(M/p^mM,g_1[m]\phi_m,\vartheta_m g_1[m]^{-1})$ and  $(M/p^mM,g_2[m]\phi_m,\vartheta_m g_2[m]^{-1})$.

\smallskip
{\bf (b)} Let $g_1,g_2\in G(W(k))$. By an {\it inner isomorphism} between the two triples $(M,g_1\phi,G)$ and  $(M,g_2\phi,G)$, we mean an element $g_3\in G(W(k))$ that is an isomorphism between the Dieudonn\'e modules  $(M,g_1\phi)$ and  $(M,g_2\phi)$ (i.e., we have $g_3g_1\phi=g_2\phi g_3$).

\smallskip
{\bf (c)} Let $g_1,g_2\in G(W(k))$. By a {\it rational inner isomorphism} between the two triples $(M,g_1\phi,G)$ and $(M,g_2\phi,G)$, we mean an element $g_3\in G(B(k))$ that is an isomorphism between $(M[{1\over p}],g_1\phi)$ and  $(M[{1\over p}],g_2\phi)$ (i.e., we have $g_3g_1\phi=g_2\phi g_3$).

\smallskip
{\bf (d)} Let $\pmb{\text{Aut}}(D[p^m])_{\text{crys},\red}$ and $\lambda_{m,\text{red}}$ be as in Theorem 2.3. Let $\pmb{\text{Aut}}(D[p^m])^G_{\text{crys},\red}$ be the reduced group of the intersection $\scrD_m^G\cap \pmb{\text{Aut}}(D[p^m])_{\text{crys},\red}$ taken inside $\scrD_m$. Let $\pmb{\text{Aut}}(D[p^m])^G_{\text{red}}$ be the unique reduced subgroup of $\pmb{\text{Aut}}(D[p^m])_{\text{red}}$ with the property that: 
$$\pmb{\text{Aut}}(D[p^m])^G_{\text{red}}(k)=\{a\in\pmb{\text{Aut}}(D[p^m])_{\text{red}}(k)|\lambda_{m,\text{red}}(k)(a)\in \pmb{\text{Aut}}(D[p^m])^G_{\text{crys},\red}(k)\}.$$
Thus $\pmb{\text{Aut}}(D[p^m])^G_{\text{red}}$ is the reduced group of $\lambda_{m,\text{red}}^*(\pmb{\text{Aut}}(D[p^m])^G_{\text{crys},\red})$ and $\pmb{\text{Aut}}(D[p^m])^G_{\text{red}}(k)$ is the subgroup of $\pmb{\text{Aut}}(D[p^m])_{\text{red}}(k)$ formed by those elements that define (via $\lambda_{m,\text{red}}$) {\it inner automorphisms} of $(M/p^mM,\phi_m,\vartheta_m,G_{W_m(k)})$. 

\smallskip
{\bf (e)} By the {\it centralizing $G$-sequence} of $D$, we mean the sequence $(\gamma_D^G(m))_{m\ge 1}$, where $\gamma_D^G(m):=\dim(\pmb{\text{Aut}}(D[p^m])^G_{\text{red}})=\dim(\pmb{\text{Aut}}(D[p^m])^G_{\text{crys},\red})$.

\smallskip
{\bf (f)} Suppose that $D$ has a polarization $\Lambda$. By the {\it centralizing sequence} of $(D,\Lambda)$, we mean the sequence $(\gamma_{D,\Lambda}(m))_{m\ge 1}$, where $\gamma_{D,\Lambda}(m):=\dim(\pmb{\text{Aut}}((D,\Lambda)[p^m])))$.

\medskip\noindent
{\bf 4.1.5. Lemma.} {\it Let $g_1,g_2\in G(W(k))$. Then the two points $g_1[m],g_2[m]\in\scrD_m^G(k)$ belong to the same orbit of the action $\dbT_m^G$ if and only if the following two quadruples $(M/p^mM,g_1[m]\phi_m,\vartheta_m g_1[m]^{-1},G_{W_m(k)})$ and  $(M/p^mM,g_2[m]\phi_m,\vartheta_m g_2[m]^{-1},G_{W_m(k)})$ are inner isomorphic.}

\medskip
\proof
As $\scrW_{+0}^G$ is the semidirect product of $\scrW_+^G$ and $\scrW_0^G$ (cf. Lemma 4.1.2) and as the product morphism $\scrP_0^G:\scrH^G\to G$ is an open embedding, each element $h\in G(W(k))$ with the property that $h$ modulo $p$ belongs to $\scrW_{+0}^G(k)$, is of the form $h=h_1h_2h_4$, where $h_1\in\scrW_+^G(W(k))$,  $h_2\in\scrW_0^G(W(k))$, and $h_4\in\Ker(\scrW_-^G(W(k))\to\scrW_-^G(k))$. Based on this, the proof of the Lemma is the same as of Lemma 2.2.1. Strictly speaking, the reference [Va2, Lem. 3.2.2] (used in the proof of Lemma 2.2.1) is stated only for the case when $u=0$ (i.e., when $G=G_1$). But the proof of loc. cit. applies entirely to our slightly more general context in which $u\in\{0,1\}$.\endproof

\medskip\smallskip\noindent
{\bf 4.2. Quasi Shimura $p$-varieties of Hodge type.} Until the end we assume that $c=d$ and we only use $d$. Let $l\ge 3$ be an integer relatively prime to $p$. Let $\scrA_{d,1,l}$ be the Mumford moduli scheme over $\dbZ_{(p)}$ that parametrizes principally polarized abelian schemes that are of relative dimension $d$ over $\dbZ_{(p)}$-schemes  and that have a level-$l$ symplectic similitude structure, cf. [MFK, Thms. 7.9 and 7.10] adapted naturally to level-$l$ symplectic similitude structures instead of only level-$l$ (symplectic) structures. Let $(\scrD_{d,1,l},\Lambda_{\scrD_{d,1,l}})$ be the  principally quasi-polarized $p$-divisible group over $\scrA_{d,1,l}$ of the universal principally polarized abelian scheme over $\scrA_{d,1,l}$.

\medskip\noindent
{\bf 4.2.1. Definition.} Suppose that $D$ has a principal quasi-polarization $\Lambda$. Let $\psi:M\times M\to W(k)$ be the perfect, alternating form on $M$ induced naturally by $\Lambda$; for $x,y\in M$ we have $\psi(\phi(x),\phi(y))=p\sigma(\psi(x,y))$. Suppose that $G$ is a closed subgroup scheme of $\pmb{Sp}(M,\psi)$. We recall that the axioms 4.1 (i) and (ii) hold for the triple $(M,\phi,G)$. As $\mu:\dbG_m\to G_1$ can not factor through $\pmb{Sp}(M,\psi)$, we have $u=1$ (i.e., we have a short exact sequence $1\to G\to G_1\to\dbG_m\to 1$). Let $\scrF:=\{(M,g\phi,\psi,G)|g\in G(W(k))\}$; it is a {\it family} of principally quasi-polarized Dieudonn\'e modules with a group over $k$. 
By a {\it quasi Shimura $p$-variety of Hodge type relative to $\scrF$}, we mean a smooth $k$-scheme $\scrM$ equipped with a quasi-finite morphism $\scrM\to {\scrA_{d,1,l}}_k$ that satisfies the following three axioms:

\medskip
{\bf (i)} the smooth $k$-scheme $\scrM$ is equidimensional of dimension $e_-$;

\smallskip
{\bf (ii)} the morphism $\scrM\to {\scrA_{d,1,l}}_k$ induces $k$-epimorphisms at the level of complete, local rings of residue field $k$ (i.e., it is a formal closed embedding at all $k$-valued points);

\smallskip
{\bf (iii)} there exists a family of \'etale maps $\rho_i:U_i\to\scrM$ indexed by a finite set $\scrI$ for which the following four properties hold:

\medskip
\item{\bf (iii.a)}  we have $\scrM=\cup_{i\in \scrI} \im(\rho_i)$ and each $U_i=\Spec(R_i)$ is an affine $k$-scheme;

\smallskip
\item{\bf (iii.b)} if $(N_i,\Phi_{N_i},\nabla_{N_i},\psi_i)$ is the projective limit indexed by positive integers $l$ of the evaluations at the thickenings $\grS_l(R_i)$ of the principally quasi-polarized $F$-crystal $\grC_i$ over $R_i$ of the pull back of $(\scrD_{d,1,l},\Lambda_{\scrD_{d,1,l}})$ to $U_i$, then there exists an isomorphism 
$$\varpi_{i}:(N_i,\Phi_{N_i},\nabla_{N_i},\psi_i)\arrowsim (M\otimes_{W(k)} W(R_i),g_i(\phi\otimes\sigma_{R_i}),\nabla_i,\psi)$$ 
for some connection $\nabla_i:M\otimes_{W(k)}  W(R_i)\to M\otimes_{W(k)} \Omega_{W(R_i)}^\wedge$ and for some element $g_i\in G(W(R_i))$ which gives birth to (i.e., lifts) a  morphism $U_i\to G_k$ whose composite with the quotient epimorphism $G_k\twoheadrightarrow G_k/\scrW_{+0k}^G$ is an \'etale morphism 
$$\eta_{i,1,-}:U_i\to G_k/\scrW_{+0k}^G=\scrD_1^G/\scrW_{+0k}^G$$ 
(here $\grS_l(R_i)$ and $\Omega_{W(R_i)}^\wedge$ are as in the beginning of Section 3);

\smallskip
\item{\bf (iii.c)} for all pairs $(i_1,i_2)\in \scrI\times \scrI$ such that the affine scheme $U_{i_1}\times_{\scrM} U_{i_2}=\Spec(R_{i_1,i_2})$ is non-empty, the isomorphism between the pulls back of $(M\otimes_{W(k)} W(R_{i_1}),g_{i_1}(\phi\otimes\sigma_{R_{i_1}}),\nabla_{i_1},\psi)$ and $(M\otimes_{W(k)} W(R_{i_2}),g_{i_2}(\phi\otimes\sigma_{R_{i_2}}),\nabla_{i_2},\psi)$ to $W(R_{i_1,i_2})$ induced naturally by $\varpi_{i_1}$ and $\varpi_{i_2}$, is defined by an element $h_{i_1,i_2}\in G(W(R_{i_1,i_2}))$;

\smallskip
\item{\bf (iii.d)} for each formally \'etale morphism $\Spec(k[[x_1,\ldots,x_{e_-}]])\to U_i$ that lifts a point $y\in U_i(k)$ and for each morphism $w_{-}^G:\Spec(W(k)[[x_1,\ldots,x_{e_-}]])\to G$ which modulo the ideal $(x_1,\ldots,x_{e_-})$ defines the identity section of $G$, which factors through the open subscheme $\im(\scrP_0^G)$ of $G$, and which defines naturally a formally \'etale morphism from $\Spec(W(k)[[x_1,\ldots,x_{e_-}]])$ to $\im(\scrP_0^G)/\scrW_{+0}^G\arrowsim\scrW_-^G$, the pull back of $\grC_i$ to $\Spec(k[[x_1,\ldots,x_{e_-}]])$ is isomorphic to 
$$(M\otimes_{W(k)} W(k)[[x_1,\ldots,x_{e_-}]],w_{-}^G(g_{y,i}\phi\otimes\Phi),\widehat{\nabla_i},\psi)$$ 
under an isomorphism that induces (via $\varpi_i$) an isomorphism 
$$(M\otimes_{W(k)} W(R_i),g_i(\phi\otimes\sigma_{R_i}),\nabla_i,\psi)\otimes_{W(R_i)} W(k[[x_1,\ldots,x_{e_-}]])\leqno (12)$$
$$\arrowsim (M\otimes_{W(k)} W(k)[[x_1,\ldots,x_{e_-}]],w_{-}^G(g_{y,i}\phi\otimes\Phi),\widehat{\nabla_i},\psi)\otimes W(k[[x_1,\ldots,x_{e_-}]])$$ 
\noindent
defined by some element of $G(W(k[[x_1,\ldots,x_{e_-}]]))$; here:

\medskip\noindent
-- $g_{y,i}\in G(W(k))$ is the pull back of $g_i\in G(W(R_i))$ via the Teichm\"uller section $\break\Spec(W(k))\hookrightarrow \Spec(W(R_i))$ defined by $y$;

\smallskip\noindent
-- $\Phi$ is the Frobenius lift of $W(k)[[x_1,\ldots,x_{e_-}]]$ that is compatible with $\sigma$ and that takes $x_i$ to $x_i^p$ for all $i\in\{1,\ldots,e_-\}$;

\smallskip\noindent
-- the tensorization with $W(k[[x_1,\ldots,x_{e_-}]])$ of the right hand side of (12) is via the $W(k)$-monomorphism $\nu:W(k)[[x_1,\ldots,x_{e_-}]]\hookrightarrow W(k[[x_1,\ldots,x_{e_-}]])$ that takes $x_i$ to $(x_i,0,\ldots)\in W(k[[x_1,\ldots,x_{e_-}]])$ for all $i\in\{1,\ldots,e_-\}$.

\medskip\noindent
{\bf 4.2.2. Simple properties.} {\bf (a)} Let $g_i[m]\in G(W_m(R_i))$ be the natural reduction of $g_i\in G(W(R_i))$. Let 
$$\eta_{i,m}:U_i\to\scrD_m^G$$ 
be the morphism defined by $g_i[m]$. Let $y\in \scrM(k)$. If $i\in \scrI$ is such that we have $\break y\in \im(U_i(k)\to\scrM(k))$, then we denote also by $y$ an arbitrary $k$-valued point of $U_i$ that maps to $y$ and we take $g_{y,i}\in G(W(k))$ to be as in the property (iii.d) of the axiom 4.2.1 (iii). Due to the property (iii.c) of the axiom 4.2.1 (iii), the inner isomorphism class of the triple $(M,g_{y,i}\phi,G)$ depends only on $y\in\scrM(k)$ and not on either $i\in \scrI$ or the choice of the point $y\in U_i(k)$ that maps to $y\in\scrM(k)$. From this and Lemma 4.1.5 we get that the orbit $\scrO_m^G(y)$ of $g_{y,i}[m]:=\eta_{i,m}(k)(y)\in\scrD_m^G$ (i.e., of the reduction modulo $p^m$ of $g_{y,i}$) depends only on $y\in\scrM(k)$. We call $(M,g_{y,i}\phi,G)$ the {\it $F$-crystal with a group attached} to the point $y\in\scrM(k)$. We also call $(M/p^mM,g_{y,i}[m]\phi_m,\vartheta_m g_{y,i}[m]^{-1},G_{W_m(k)})$ the {\it $D$-truncation modulo $p^m$} of $(M,g_{y,i}\phi,G)$ (to be compared with [Va3, Subsubsect. 1.1.1], where the case $m=1$ is considered but for the triple $(M,g_{y,i}\phi,G_1)$ instead of for $(M,g_{y,i}\phi,G)$). 

\smallskip
{\bf (b)} Each \'etale scheme over $\scrM$ is itself a quasi Shimura $p$-variety of Hodge type relative to $\scrF$. Thus locally in the \'etale topology of $\scrM$, for the (locally \'etale) study of $\scrM$ one can assume that $\scrI$ has only one element and that $\scrM=U_i=\Spec(R_i)$ is affine. One can also assume that there exists an \'etale $k$-monomorphism $k[x_1,\ldots,x_{e_-}]\hookrightarrow  R_i$ such that the ideal $(x_1,\ldots,x_{e_-})$ is mapped to the maximal ideal of $R_i$ that defines an a priori fixed point $y\in\scrM(y)=U_i(y)$.

\smallskip
{\bf (c)} We recall the well known argument that the reduction $\widehat{\nabla_{i,m}}$ modulo $p^m$ of the connection $\widehat{\nabla_i}$ on $M\otimes_{W(k)} W(k)[[x_1,\ldots,x_{e_-}]]$ is uniquely determined by the natural reductions $w_{-,m}^G\in G(W_m(k)[[x_1,\ldots,x_{e_-}]])$ and $g_{y,i}[m]\in G(W_m(k))$ of $w_{-}^G$ and $g_{y,i}$ (respectively). Let $\Phi_m$ be the reduction modulo $p^m$ of $\Phi$. The connection $\widehat{\nabla_{i,m}}$ on $M/p^mM\otimes_{W_m(k)} W_m[[x_1,\ldots,x_{e_-}]]$ satisfies the following equation 
$$X\circ w_{-,m}^G(g_{y,i}[m]\phi_m\otimes\Phi_m)=(w_{-,m}^G(g_{y,i}[m]\phi_m\otimes\Phi_m)\otimes d\Phi_m)\circ X$$ 
between maps from $(M+{1\over p}F^1)/p^m(M+{1\over p}F^1)\otimes_{W_m(k)} W_m[[x_1,\ldots,x_{e_-}]]$ to $M/p^mM\otimes_{W_m(k)} \oplus_{j=1}^{e_-} W_m[[x_1,\ldots,x_{e_-}]]dx_j$. Let $X_1,X_2$ be two solutions of this equation in $X$. We have $X_1-X_2\in \End(M/p^mM)\otimes_{W_m(k)}\oplus_{j=1}^{e_-}  W_m(k)[[x_1,\ldots,x_{e_-}]]dx_j$. As for all $j\in\{1,\ldots,e_-\}$ we have $d\Phi(x_j)=px_j^{p-1}dx_j$ and as $\phi_m$ maps $(M+{1\over p}F^1)/p^m(M+{1\over p}F^1)$ isomorphically to $M/p^mM$, by induction on $q\in\dbN$ we get that $X_1-X_2$ belongs to $\End(M/p^mM)\otimes_{W_m(k)}\oplus_{j=1}^{e_-}  (x_1,\ldots,x_{e_-})^qdx_j$. As the local ring $W_m(k)[[x_1,\ldots,x_{e_-}]]$ is complete in the $(x_1,\ldots,x_{e_-})$-adic topology, we conclude that $X_1-X_2=0$. Thus indeed the connection $\widehat{\nabla_{i,m}}$ is uniquely determined by $w_{-,m}^G$ and $g_{y,i}[m]$.

\medskip\noindent
{\bf 4.2.3. Level $m$ stratification of $\scrM$.} Let $y\in\scrM(k)$. Let $(y_1,y_2)\in \scrI\times \scrI$ be such that the affine scheme $U_{i_1}\times_{\scrM} U_{i_2}$ is non-empty. Due to the property (iii.c) of the axiom 4.2.1 (iii), the reduced schemes of the pulls back of $\scrO_m^G(y)$ via the two morphisms $U_{i_1}\times_{\scrM} U_{i_2}\to\scrD_m$ defined naturally by the morphisms $\eta_{i_1,m}$ and $\eta_{i_2,m}$, are equal. This implies that there exists a unique reduced, locally closed subscheme $\grs_y^G(m)$ of $\scrM$ such that we have:
\medskip
{\bf (*)} for all $i\in \scrI$, the scheme $U_i\times_{\scrM} \grs_y^G(m)$ is the reduced scheme of $\eta_{i,m}^*(\scrO_m^G(y))$. 

\medskip
Thus we have an identity of sets
$$\grs_y^G(m)(k)=\{z\in\scrM(k)|\scrO_m^G(z)=\scrO_m^G(y)\}.$$ 
For each algebraically closed field $k_1$ that contains $k$ and for each point $y_1\in\scrM(k_1)$, we similarly define an orbit $\scrO_m^G(y_1)$ of the action $\dbT_m^{G_{W(k_1)}}=\dbT_m^G\times_k k_1$ and a reduced, locally closed subscheme $\grs_{y_1}^G(m)$ of $\scrM_{k_1}$ such that we have an identity $\grs_{y_1}^G(m)(k_1)=\{z_1\in\scrM(k_1)|\scrO_m^G(z_1)=\scrO_m^G(y_1)\}$ of sets of $k_1$-valued points.

By the {\it level $m$ stratification} $\scrS^G(m)$ of $\scrM$, we mean the stratification in the sense of [Va2, Def. 2.1.1] defined by the rule:
\medskip
{\bf (**)} for each algebraically closed field $k_1$ that contains $k$, the set of strata of $\scrS^G(m)$ that are reduced, locally closed subschemes of $\scrM_{k_1}$ is $\scrS^G_{k_1}(m):=\{\grs_{y_1}^G(m)|y_1\in\scrM(k_1)\}$.

\medskip\noindent
{\bf 4.2.4. Ultimate stratification of $\scrM$.}
Let $n_G$ be the smallest positive integer that has the following property (cf. [Va2, Main Thm. A]):

\medskip
{\bf (*)} for each $g\in G(W(k))$ and every $g_{n_G}\in\Ker(G(W(k))\to G(W_{n_G}(k)))$, there exists an inner isomorphism between $(M,g\phi,G)$ and  $(M,g_{n_G}g\phi,G)$.

\medskip
Based on (*), for all $m\ge n_G$ we have $\scrS^G(m)=\scrS^G(n_G)$. We call 
$$\scrS^G:=\scrS^G(n_G)$$ 
the {\it ultimate stratification} or the {\it Traverso stratification} of $\scrM$.

\medskip\noindent
{\bf 4.3. Basic Corollary.} {\it Let $(M,\phi,\psi)$ be the principally quasi-polarized Dieudonn\'e module of a principally quasi-polarized $p$-divisible group $(D,\Lambda)$ over $k$ of height $2d$. Let $G$ be a smooth, closed subgroup scheme of $\pmb{Sp}(M,\psi)$ such that the axioms 4.1 (i) and (ii) hold for the triple $(M,\phi,G)$ (with $u=1$). Let $l\ge 3$ be a positive integer prime to $p$. Let $\scrM\to {\scrA_{d,1,l}}_k$ be a quasi Shimura $p$-variety of Hodge type relative to $\scrF:=\{(M,g\phi,\psi,G)|g\in G(W(k))\}$. To fix the notations, we assume that there exists a point $y\in\scrM(k)$ such that (with the notations of Subsubsections 4.1.3 and 4.2.1 (a)) we have $\scrO_m^G(y)=\scrO_m^G$. Then the following five properties hold:

\medskip
{\bf (a)} the reduced, locally closed subscheme $\grs_y^G(m)$ of $\scrM$ is regular and equidimensional;

\smallskip
{\bf (b)} we have $\dim(\grs_y^G(m))=e_--\gamma_D^G(m)$;

\smallskip
{\bf (c)} for $m\ge n_G$, the number $\gamma_D^G(m)$ is equal to $\gamma_D^G(n_G)$; 

\smallskip
{\bf (d)} the $\scrM$-scheme $\grs_y^G(m)$ is quasi-affine;

\smallskip
{\bf (e)} if the image of the abstract composite homomorphism 
$$\chi_D^G(m):\pmb{\text{Aut}}(D[p^m])^G(k)\to\pmb{\text{Aut}}(D[p])^G(k)\to\scrW_{+0k}^G(k)$$
$$\twoheadrightarrow (\scrW_{+0k}^G/\scrW_{+k}^G\scrW_{0k}^{G,\text{unip}})(k)\arrowsim (\scrW_{0k}^G/\scrW_{0k}^{G,\text{unip}})(k)$$ is finite, then the reduced, locally closed subscheme $\grs_y^G(m)$ of $\scrM$ satisfies the {\it purity property} (here the second homomorphism $\pmb{\text{Aut}}(D[p])^G(k)\to\scrW_{+0k}^G(k)$ is defined naturally by the homomorphism $\lambda_{1,\text{red}}$ of Theorem 2.4 (b) and $\scrW_{0k}^{G,\text{unip}}$ is the unipotent radical of $\scrW_{0k}^G$).}

\medskip
\proof
We check (a). Let $y_1,y_2\in\grs_y^G(m)(k)$. For $j\in\{1,2\}$, let $I_j^G$ be the completion of the local ring of $\scrM$ at $y_j$ and let $i_j\in \scrI$ be such that there exists a point (denoted in the same way) $y_j\in U_{i_j}(k)$ that maps to $y_j\in\scrM(k)$. Due to the axiom 4.2.1 (i), we can identify $I_1^G=I_2^G=k[[x_1,\ldots,x_{e_-}]]$. Let $\grs_j^G:=\Spec(I^G_j)\times_{\scrM} \grs_y^G(m)$; it is a reduced, local, complete, closed subscheme of $\Spec(I_j^G)$. As $y_1,y_2\in\grs_y^G(m)(k)$, from Lemma 4.1.5 we get that (up to inner isomorphisms) we can assume that $g_{y_1,i_1}[m]=g_{y_2,i_2}[m]$. This assumption is compatible with the property (iii.d) of the axiom 4.2.1 (iii); this is so as the conjugate of each morphism $w_-^G:\Spec(W(k)[[x_1,\ldots,x_{e_-}]])\to G$ as in the mentioned property via an element of $\im(\scrP_{0-}^G(W(k)))=\{\flat\in G(W(k))|\flat [1]\in\scrW_{+0}^G(k)\}$, is a morphism that has the same properties as $w_-^G$. From this and the property (iii.d) of the axiom 4.2.1 (iii) we get that there exists an isomorphism $\gamma_{12}:\Spec(I_1^G)\arrowsim\Spec(I_2^G)$ such that the reduction modulo $p^m$ of the pull back of $(M\otimes_{W(k)} W(k)[[x_1,\ldots,x_{e_-}]],w_{-}^G(g_{y_2,i_2}\phi\otimes\Phi_{W(k)[[x_1,\ldots,x_{e_-}]]}),\widehat{\nabla_{i_2}},\psi)$ via $\gamma_{12}$ is naturally identified with the reduction modulo $p^m$ of $(M\otimes_{W(k)} W(k)[[x_1,\ldots,x_{e_-}]],w_{-}^G(g_{y_1,i_1}\phi\otimes\Phi_{W(k)[[x_1,\ldots,x_{e_-}]]}),\widehat{\nabla_{i_1}},\psi)$ (cf. also Subsubsection 4.2.2 (c) for the part involving connections). Such a natural identification $\scrI_{12}$ is unique up to isomorphisms defined by elements of $G(W_m[[x_1,\ldots,x_{e_-}]])$. 

Let $K$ be an algebraically closed field that contains $k$. Let $z_1\in\Spec(I_1^G)(K)$ and let $z_2:=\gamma_{12}\circ z_1\in\Spec(I_2^G)(K)$. Let $g_{z_j}\in G(W(K))$ be such that the triple $(M\otimes_{W(k)} W(K),g_{z_j}(\phi\otimes\sigma_K),G_{W(K)})$ is the $F$-crystal with a group over $K$ attached to the $K$-valued point of $\scrM$ defined naturally by $z_j$. Let $\Theta_j:\Spec(W(K))\to \Spec(W(k)[[x_1,\ldots,x_{e_-}]])$ be the Teichm\"uller lift of the point $z_j\in\Spec(I_j^G)(K)=\Spec(k[[x_1,\ldots,x_{e_-}]])(K)$. The element $g_{z_j}\in G(W(K))$ is the composite of $\Theta_j$ with $w_{-}^Gg_{y_j,i_j}\in G(W(k)[[x_1,\ldots,x_{e_-}]])$. Due to the existence of the identification $\scrI_{12}$, we get that we can assume that $g_{z_1}$, $g_{z_2}\in G(W(K))$ are congruent modulo $p^m$. This implies that $z_1$ factors through $\grs_1^G$ if and only if $z_2$ factors through $\grs_2^G$. Thus we have an identity $\grs_1^G=\gamma_{12}^*(\grs_2^G)$. 

As in Subsection 3.2 we argue that the identity $\grs_1^G=\gamma_{12}^*(\grs_2^G)$ implies that $\grs_y^G(m)$ is a regular and equidimensional $k$-scheme. Thus (a) holds. 

We check (b). The finite epimorphism $\iota_m:\scrC_m\twoheadrightarrow\pmb{\text{Aut}}(D[p^m])_{\text{crys},\text{red}}$ induces a finite homomorphism $\iota_m^G:\scrC_m^G\to\pmb{\text{Aut}}(D[p^m])_{\text{crys},\text{red}}^G$, cf. the very definitions of $\scrC_m^G$ and $\pmb{\text{Aut}}(D[p^m])^G_{\text{crys},\text{red}}$. The homomorphism $\iota_m^G(k):\scrC_m^G(k)\to\pmb{\text{Aut}}(D[p^m])_{\text{crys},\text{red}}^G(k)$ is injective (cf. Theorem 2.4 (b)) and next we will check that it is also surjective. Let $h_1[m]h_2[m]h_3[m]^p\in\pmb{\text{Aut}}(D[p^m])_{\text{crys},\text{red}}^G(k)\leqslant G(W_m(k))$, where $(h_1[m],h_2[m],h_3[m])\in\scrC_m(k)$ (cf. Theorem 2.4 (b)). As $h_1[m]h_2[m]h_3[m]^p\in G(W_m(k))$, we have the relations $h_1[m]\in\scrW_+^G(W_m(k))$, $h_2[m]\in\scrW_0^G(W_m(k))$, and $h_3[m]^p\in\scrW_-^G(W_m(k))$. But as $(h_1[m],h_2[m],h_3[m])\in\scrC_m(k)$, we have $h_1[m]h_2[m]h_3[m]^p=\sigma_{\phi}(h_1[m]^p)\sigma_{\phi}(h_2[m])\sigma_{\phi}(h_3[m])$. From the last two sentences we get that $\sigma_{\phi}(h_3[m])\in G(W_m(k))$ and therefore (cf. Subsubsection 4.1.3) that $h_3[m]\in G(W_m(k))$. Thus $h_3[m]\in\scrW_-^G(W_m(k))=G(W_m(k))\cap \scrW_-(W_m(k))$ and therefore $(h_1[m],h_2[m],h_3[m])\in\scrC_m^G(k)$. Thus the homomorphism $\iota_m^G(k)$ is an isomorphism. Therefore we have the following analogue 
$$\dim(\scrC_m^G)=\dim(\pmb{\text{Aut}}(D[p^m])_{\text{crys},\text{red}})=\gamma^G_D(m)$$
of Theorem 2.4 (c) (for the last identity cf. Definitions 4.1.4 (d) and (e)). 

As $\dim(\scrC_m^G)=\gamma^G_D(m)$ and as for each $i\in \scrI$ the morphism $\eta_{i,1,-}:U_i\to \scrD_1^G/\scrW_{+0k}^G$ is \'etale (cf. property (iii.b) of the axiom 4.2.1 (iii)), the proof of (b) is the same as the proof of Theorem 1.2 (c) presented in Subsection 3.3.  

We check (c). For $m\ge n_G$, we have $\grs_y^G(m)=\grs_y^G(n_G)$ (cf. the definition of $n_G$). From this and (b) we get that we have $\gamma_D^G(m)=\gamma_G^G(n_G)$. 

We check (d). As $\scrO_m^G$ is a quasi-affine $\scrD_m^G$-scheme, for each $i\in \scrI$ the scheme $\grs_y^G(m)\times_{\scrM} U_i$ is a quasi-affine $U_i$-scheme. From this and the \'etaleness part of the axiom 4.2.1 (iii) we get that  $\grs_y^G(m)$ is a quasi-affine $\scrM$-scheme. Thus (d) holds. 

We check (e). From the hypotheses of (e) and from Definition 4.1.4 (d), as in Lemma 2.5 we argue that $\scrC_m^{0G}$ is a subgroup of $\Ker(\Xi_m^G)$, where the epimorphism 
$$\Xi_m^G:\scrH_m^G\twoheadrightarrow \scrW_{0k}^G/\scrW_{0k}^{G,\text{unip}}$$ 
is defined at the level of $k$-valued points by the rule: $(h_1[m],h_2[m],h_3[m])\in\scrH_m^G(k)$ is mapped to the image of $h_2[1]\in\scrW_{0k}^G(k)$ in $(\scrW_{0k}^G/\scrW_{0k}^{G,\text{unip}})(k)$. From this and Lemma 2.2.3 we get that $\scrC_m^{0G}$ is a subgroup of the unipotent radical of $\scrH_m^G$. Thus the proof of Lemma 2.6 adapts to give us that the orbit $\scrO_m^G\arrowsim\scrH_m^G/\scrS_m^G$ is an affine scheme. This implies that for each $i\in \scrI$ the scheme $\grs_y^G(m)\times_{\scrM} U_i$ is an affine $U_i$-scheme. From this and the \'etaleness part of the axiom 4.2.1 (iii) we get that  $\grs_y^G(m)$ is an affine $\scrM$-scheme. Thus (e) holds. 
\endproof

\medskip\noindent
{\bf 4.3.1. Corollary.} {\it {\bf (a)} If there exists a point $y\in\scrM(k)$ such that we have $\scrO_m^G(y)=\scrO_m^G$, then $\gamma_D^G(m)\le e_-$.

\smallskip
{\bf (b)} The level $n_G$ stratification $\scrS^G(n_G)$ of $\scrM$ satisfies the purity property as defined in [Va2, Subsubsect. 2.1.1].}

\medskip
\proof
Part (a) follows from Corollary 4.3 (b). For $m>>n_G$, the images of $\pmb{\text{Aut}}(D[p^m])(k)$ and $\text{Aut}(D)$ in $\pmb{\text{Aut}}(D[p])(k)$ are equal (cf. [Va2, Thm. 5.1.1 (c)]) and thus they are finite. Therefore the hypotheses of Corollary 4.3 (e) hold for $m>>n_G$. Thus $\grs_y^G(m)=\grs_y^G(n_G)$ is an affine $\scrM$-scheme, cf. Corollary 4.3 (e). Part (b) is only a reformulation of the last sentence in terms of stratifications.
\endproof

\medskip\noindent
{\bf 4.3.2. Remarks.} {\bf (a)} Suppose that there exits $q\in\dbN$ such that the ${\scrA_{d,1,l}}_k$-scheme $\scrM$ is the pull back of an ${\scrA_{d,1,l}}_{\dbF_{p^q}}$-scheme $\scrM_{\dbF_{p^q}}$. If the axiom 4.2.1 (iii) holds naturally over ${\dbF_{p^q}}$, then the stratification $\scrS^G(m)$ is the pull back of a stratification of $\scrM_{\dbF_{p^q}}$. 

\smallskip
{\bf (b)} The proof of Theorem 1.2 (e) can be easily adapted to show that the number $\gamma_D^G(n_G)$ depends only on the rational inner isomorphism class of $(M,\phi,G)$. 

\medskip\smallskip\noindent
{\bf 4.4. Example.} Suppose that $\scrW_{0k}^G/\scrW_{0k}^{G,\text{unip}}$ is a torus over $k$. As $\scrC_m^0$ is a unipotent group (cf. Theorem 2.4 (a)), the image of the subgroup $\scrC_m^{0G}$ of $\scrC_m^0$ in $\scrW_{0k}^G/\scrW_{0k}^{G,\text{unip}}$ is trivial. Thus the hypothesis of Corollary 4.3 (e) holds for all $m\ge 1$ and therefore the stratification $\scrS^G(m)$ satisfies the purity property. On the other hand, we emphasize that the number $n_G$ can be arbitrarily large (to be compared with [Va2, Subsect. 4.4]). 

The previous paragraph applies if $G_1$ is a reductive group scheme whose adjoint group scheme $G_1^{\ad}$ is isomorphic to $\pmb{PGL}_2^t$ for some $t\in\dbN$ and if the cocharacter $\mu:\dbG_m\to G_1$ has a non-trivial image in each simple factor $\pmb{PGL}_2$ of $G_1^{\ad}$; in such a case $\scrW_0^G$ is itself a torus of rank $t$ over $\Spec(W(k))$ and therefore the group $\scrW_{0k}^G/\scrW_{0k}^{G,\text{unip}}$ is a torus over $k$.  

\medskip\smallskip\noindent
{\bf 4.5. Example.} Suppose that $G=\pmb{Sp}(M,\psi)$ and $G_1=\pmb{GSp}(M,\psi)$. Let $M=F^1\oplus F^0$ be a direct sum decomposition  such that $F^1/pF^1$ is the kernel of $\phi$ modulo $p$ and we have $\psi(F^1,F^1)=\psi(F^0,F^0)=0$. Obviously both axioms 4.1 (i) and (ii) hold for the triple $(M,\phi,G)$. Let $\scrM:={\scrA_{d,1,l}}_k$. It is well known that $d_G=\dim(G_{B(k)})=2d^2+d$, that $e_0=d^2$, that $e_-=e_+={{d(d+1)}\over 2}$, and that $\dim(\scrM)$ is a smooth $k$-scheme which is equidimensional of dimension $e_-$. Thus the axioms 4.2.1 (i) and (ii) hold. It is easy to check that the axiom 4.2.1 (iii) holds as well. For instance, we can take $(U_i)_{i\in \scrI}$ to be an arbitrary affine, open cover of $\scrM$ such that for all $i\in \scrI$ the following three properties hold: (i) the $R_i$-module $N_i/pN_i$ is free of rank $2d$, (ii) the kernel of the reduction of $\Phi_{N_i}$ modulo $p$ is a free direct summand of $N_i/pN_i$ of rank $d$, and (iii) we have an \'etale morphism $k[x_1,\ldots,x_{e_-}]\hookrightarrow R_i$; here $N_i$, $\Phi_{N_i}$, and $R_i$ are as in the property (iii.b) of the axiom 4.2.1 (iii). The fact that this property holds is an easy consequence of the following fact: (iv) any two symplectic spaces over $R_i$ that are defined by free $R_i$-modules of rank $2d$ and that have free direct summands which are of rank $d$ and isotropic, are isomorphic; thus as in the proof of the property 3.1 (v) we argue that (by shrinking $U_i$) we can assume that the morphism $\eta_{i,1,-}:U_i\to\scrD_1^G/\scrW_{+0k}^G$ is \'etale. The fact that the property (iii.d) of the axiom 4.2.1 (iii) holds, is a particular case of the Faltings deformation theory in the form used in [Va1, Subsect. 5.4] (loc. cit. worked with $p>2$ but it also applies even if $p=2$). 

Thus $\scrM$ is a quasi Shimura $p$-variety of Hodge type relative to $\scrF$. The level $m$ stratification $\scrS^G(m)$ of $\scrM$ is the pull back of a stratification of ${\scrA_{d,1,l}}_{\dbF_p}$, cf.  Remark 4.3.2 (a).

Each principally quasi-polarized truncated Barsotti--Tate group $(B,\Lambda_B)$ of level $m$ and dimension $d$ over $k$ lifts to a principally quasi-polarized $p$-divisible group over $k$. From this and [Va2, Prop. 5.3.3] we get that $(B,\Lambda_B)$ is the principally quasi-polarized truncated Barsotti--Tate group of level $m$ of a principally polarized abelian variety of dimension $d$ over $k$. 
Each principally quasi-polarized Dieudonn\'e module of dimension $d$ over $k$ is isomorphic to $(M,u\phi,\psi)$ for some element $u\in G(W(k))$. Based on the last two sentences and on Lemma 4.1.5, we get that the strata of $\scrS^G_k(m)$ (equivalently, the orbits of the action $\dbT_m^G$) are parametrized by {\it isomorphism classes of principally quasi-polarized truncated Barsotti--Tate groups of level $m$ and dimension $d$ over $k$}. 

To fix the notations, let $y\in\scrM(k)$ be such that $y^*((\scrD_{d,1,l},\Lambda_{\scrD_{d,1,l}})_{\scrM})=(D,\Lambda)$. Let $n_{D,\Lambda}$ be the smallest positive integer such that $(D,\Lambda)$ is uniquely determined up to isomorphisms by its truncation $(D,\Lambda)[p^{n_{D,\Lambda}}]$ of level $n_{D,\Lambda}$, cf. [Va2, Subsubsect. 3.2.5]. For all $m\ge 1$ we have $\pmb{\text{Aut}}(D[p^m])^G_{\text{red}}=\pmb{\text{Aut}}((D,\Lambda)[p^m])_{\text{red}}$, cf. the classical Dieudonn\'e theory. Thus for all $m\ge 1$ we have $\gamma_D^G(m)=\gamma_{D,\Lambda}(m)$. If $m\ge n_{D,\Lambda}$, then $\grs_y^G(m)=\grs_y^G(n_{D,\Lambda})$. The proof of Theorem 1.2 (f) can be easily adapted to provide a formula for $\gamma_{D,\Lambda}(n_{D,\Lambda})$ in terms of the Newton polygons slopes of $(M,\phi)$. As this formula for $\gamma_{D,\Lambda}(n_{D,\Lambda})$ is stated in an informal manuscript of Oort, it will not be presented here. 

\medskip\smallskip\noindent
{\bf 4.6. Example.} For basic properties of Shimura varieties of Hodge type we refer to [De1], [De2], [Mi1], [Mi2], and [Va1]. Let $(W,\psi)$ be a symplectic space over $\dbQ$ of dimension $2d$. Let $\dbS$ be the $2$-dimensional torus over $\dbR$ whose group of $\dbR$-valued points is $\dbG_m(\dbC)$. Let $\scrY$ be the set of all homomorphisms $\dbS\to \pmb{GSp}(W,\psi)_{\dbR}$ that define Hodge $\dbQ$--structures on $W$ of type $\{(-1,0),(0,-1)\}$ and that have either $2\pi i\psi$ or $-2\pi i\psi$ as polarizations. Let $(\scrG,\scrX)\hookrightarrow (\pmb{GSp}(W,\psi),\scrY)$ be an injective map of Shimura pairs. Let $\Sh(\scrG,\scrX)$ and $\Sh(\pmb{GSp}(W,\psi),\scrY)$ be the canonical models of $(\scrG,\scrX)$ and $(\pmb{GSp}(W,\psi),\scrY)$ over the reflex fields $E(\scrG,\scrX)$ and $\dbQ$ (respectively), cf. [De1, Variant 5.9]. Let $w$ be a prime of $E(\scrG,\scrX)$ that divides $p$. Let $O_{(w)}$ be the localization of the ring of integers of $E(\scrG,\scrX)$ with respect to $w$. Let $L$ be a $\dbZ$-lattice of $W$ such that we have a perfect, alternating form $\psi:L\times L\to\dbZ$. 

Let $l\ge 3$ be an integer prime to $p$. Let $K(l):=\{h\in \pmb{GSp}(W,\psi)(\widehat{\dbZ})|h\,\,\text{ modulo}\,\,l\,\,\text{is}\,\,\text{the}\,\,\text{identity}\}$. Let $H(l)$ be an open subgroup of $\scrG(\widehat{\dbZ}\otimes_{\dbZ} \dbQ)\cap K(l)$. There exists a natural finite morphism
$$\Sh(\scrG,\scrX)/H(l)\to\Sh(\pmb{GSp}(W,\psi),\scrY)_{E(\scrG,\scrX)}/K(l)$$ 
\noindent
(cf. [De1, Cor. 5.4 and Def. 3.13]). We also have a natural identification $\Sh(\pmb{GSp}(W,\psi),\scrY)/K(l)={\scrA_{d,1,l}}_{\dbQ}$, cf. [De1, Example 4.16]. Thus we can speak about the normalization $\scrN$ of ${\scrA_{d,1,l}}_{O_{(w)}}$ in the ring of fractions of $\Sh(\scrG,\scrX)/H(l)$. Let $\scrM$ be a connected component of $\scrN_k$. We assume that the following two properties hold:

\medskip
{\bf (i)} the prime $w$ is unramified over $p$ and the $O_{(w)}$-scheme $\scrN$ is smooth; 

\smallskip
{\bf (ii)} the natural morphism $\scrM\to {\scrA_{d,1,l}}_k$ induces $k$-epimorphisms at the level of complete, local rings of residue field $k$ (i.e., it is a formally closed embedding at all $k$-valued points).

\medskip
One can check that (i) implies (ii) (to be compared with [Va1, Cor. 5.6.1]). Due to the first part of the property (i), we can view $W(k)$ as an $O_{(w)}$-algebra. Let $y\in\scrM(k)$. Let $z\in\scrN(W(k))=\scrN_{W(k)}(W(k))$ be a point that lifts $y$. Let $(M,F^1,\phi,\psi)$ be the principally quasi-polarized filtered Dieudonn\'e module over $k$ of $z^*((\scrD_{d,1,l},\Lambda_{\scrD_{d,1,l}})_{\scrN})$. Let $G_{1B(k)}$ be the connected subgroup of $\pmb{GL}_{M[{1\over p}]}$ that corresponds naturally to $\scrG$ via Fontaine comparison theory, as in [Va1, Subsubsects. 5.3.4 and 5.6.5]. Let $G_1$ be the Zariski closure of $G_{B(k)}$ in $\pmb{GL}_M$. Let $G:=G_1\cap \pmb{Sp}(M,\psi)$. We also assume that the following third property holds:

\medskip
{\bf (iii)} the group scheme $G_1$ is smooth over $\Spec(W(k))$.

\medskip
It is easy to see that the axioms 4.1 (i) and (ii) hold for the triple $(M,\phi,G)$ and that $e_-=\dim(\scrX)=\dim(\scrM)$ (to be compared with [Va1, Subsubsects. 5.4.6 and 5.4.7]). We have a short exact sequence $1\to G\to G_1\to\dbG_m\to 1$ which, due to the existence of a cocharacter $\mu:\dbG_m\to G_1$ as in the axiom 4.1 (ii),  splits. This implies that $G_1$ is the semidirect product of $G$ and $\dbG_m$. Based on this and the property (iii), we get that $G$ is smooth over $\Spec(W(k))$. One can use Faltings deformation theory as in [Va1, Subsect. 5.4], to check that $\scrM$ is a quasi Shimura $p$-variety of Hodge type relative to $\scrF:=\{(M,g\phi,\psi,G)|g\in G(W(k))\}$ and that, up to natural identifications, $\scrF$ does not depend on the choice of the point $y\in\scrM(k)$. The \'etale morphisms $U_i\to\scrM$ are obtained via Faltings deformation theory through Artin's approximation theory.  Thus $\scrM$ has a level $m$ stratification $\scrS^G(m)$ which has all the properties described in Basic Corollary 4.3.

\medskip\noindent
{\bf 4.6.1. Example.} If $\Sh(\scrG,\scrX)$ is a {\it Hilbert--Blumenthal moduli variety} (i.e., if the adjoint group $\scrG^{\ad}_{\dbR}$ is isomorphic to $\pmb{PGU}(1,1)^t$ for some $t\in\dbN$) and if $G_1$ is a reductive group scheme over $\Spec(W(k))$, then $G_1^{\ad}$ is isomorphic to $\pmb{\text{PGL}}_2^t$ and every cocharacter $\mu:\dbG_m\to G_1$ as in the axiom 4.1 (ii) has a non-trivial image in each simple factor of $G_1^{\ad}$. Thus for all $m\ge 1$ the stratification $\scrS^G(m)$ satisfies the purity property, cf. Example 4.4. 

\medskip\noindent
{\bf Acknowledgments.} We would like to thank University of Arizona, University of Bielefeld, and Max-Planck Institute, Bonn for good conditions with which to write this paper. We would also like to thank J. E. Humphreys for mentioning to us the reference [CPS] and the referee for many valuable comments and suggestions. 

\bigskip
\references{37}
{\nspace{

\Ref[Be] P. Berthelot, 
\sl Cohomologie cristalline des sch\'emas de caract\'eristique $p>0$, 
\rm Lecture Notes in Math., Vol. {\bf 407}, Springer-Verlag, Berlin-New York, 1974.

\Ref[Bo]
A. Borel,
\sl Linear algebraic groups. Second edition,
\rm Grad. Texts in Math., Vol. {\bf 126}, Springer-Verlag, New York, 1991.

\Ref[BBM] 
P. Berthelot, L. Breen, and W. Messing, 
\sl Th\'eorie de Dieudonn\'e crystalline II, 
\rm Lecture Notes in Math., Vol. {\bf 930}, Springer-Verlag, Berlin, 1982.

\Ref[BLR]
S. Bosch, W. L\"utkebohmert, and M. Raynaud, 
\sl N\'eron models, 
\rm Ergebnisse der Mathematik und ihrer Grenzgebiete (3), Vol. {\bf 21}, Springer-Verlag, Berlin, 1990.

\Ref[CPS]
E. Cline, B. Parshall, and L. Scott,
\sl Induced modules and affine quotients,
\rm Math. Annalen {\bf 230} (1977), no. 1, pp. 1--14.

\Ref[De1]
P. Deligne,
\sl Travaux de Shimura,
\rm S\'eminaire  Bourbaki, 23\`eme ann\'ee (1970/71), Exp. No. 389, pp. 123--165, Lecture Notes in Math., Vol. {\bf 244}, Springer-Verlag, Berlin, 1971.

\Ref[De2]
P. Deligne,
\sl Vari\'et\'es de Shimura: interpr\'etation modulaire, et
techniques de construction de mod\`eles canoniques,
\rm Automorphic forms, representations and $L$-functions (Oregon State Univ., Corvallis, OR, 1977), Part 2,  pp. 247--289, Proc. Sympos. Pure Math., {\bf 33}, Amer. Math. Soc., Providence, RI, 1979.

\Ref[Dem]
M. Demazure, 
\sl Lectures on $p$-divisible groups, 
\rm Lecture Notes in Math., Vol. {\bf 302}, Springer-Verlag, Berlin-New York, 1972.

\Ref[Di] 
J. Dieudonn\'e, 
\sl Groupes de Lie et hyperalg\`ebres de Lie sur un corps de caract\'erisque $p>0$ (VII), 
\rm Math. Annalen {\bf 134} (1957), no. 2, pp. 114--133.

\Ref[DG]
M. Demazure, A. Grothendieck, et al. 
\sl Sch\'emas en groupes. Vol. {\bf II--III}, 
\rm S\'eminaire de G\'eom\'etrie Alg\'ebrique du Bois Marie 1962/64 (SGA 3), Lecture Notes in Math., Vol. {\bf 152--153}, Springer-Verlag, Berlin-New York, 1970. 

\Ref[Fo]
J.-M. Fontaine, 
\sl Groupes $p$-divisibles sur les corps locaux, 
\rm J. Ast\'erisque {\bf 47/48}, Soc. Math. de France, Paris, 1977.

\Ref[Gr]
M. Greenberg,
\sl Schemata over local rings,
\rm Ann. of Math. {\bf 73} (1961), no. 3, pp. 624--648.

\Ref[Gro]
A. Grothendieck, 
\sl \'El\'ements de g\'eom\'etrie alg\'ebrique. II. \'Etude globale \'el\'ementaire de quelques classes de morphisms, 
\rm Inst. Hautes \'Etudes Sci. Publ. Math., Vol. {\bf 8}, 1961.

\Ref[Il]
L. Illusie, 
\sl D\'eformations des groupes de Barsotti--Tate (d'apr\`es A. Grothendieck), 
\rm Seminar on arithmetic bundles: the Mordell conjecture (Paris, 1983/84), pp. 151--198, J. Ast\'erisque {\bf 127}, Soc. Math. de France, Paris, 1985.

\Ref[Ko]
R. E. Kottwitz,
\sl Points on some Shimura Varieties over finite fields,
\rm J. of Amer. Math. Soc. {\bf 5} (1992), no. 2, pp. 373--444.

\Ref[Kr] H. Kraft, 
\sl Kommutative algebraische p-Gruppen (mit Anwendungen auf p-divisible Gruppen und abelsche Variet\"aten), 
\rm manuscript 86 pages, Univ. Bonn, 1975.

\Ref[Ma]
J. I. Manin, 
\sl The theory of formal commutative groups in finite characteristic, 
\rm Russian Math. Surv. {\bf 18} (1963), no. 6, pp. 1--83.

\Ref[Mi1]
J. S. Milne,
\sl The points on a Shimura variety modulo a prime of good
reduction,
\rm The Zeta function of Picard modular surfaces, pp. 151--253, Univ. Montr\'eal Press, Montreal, Quebec, 1992.

\Ref[Mi2]
J. S. Milne,
\sl Shimura varieties and motives,
\rm Motives (Seattle, WA, 1991), pp. 447--523, Proc. Sympos. Pure Math., Vol. {\bf 55}, Part 2, Amer. Math. Soc., Providence, RI, 1994.

\Ref[Mo]
B. Moonen,
\sl Group schemes with additional structures and Weyl group cosets,
\rm Moduli of abelian varieties (Texel Island, 1999), pp. 255--298, Progr. of Math., Vol. {\bf 195}, Birkh\"auser, Basel, 2001.

\Ref[MFK]
D. Mumford, J. Fogarty, and F. Kirwan, 
\sl Geometric invariant theory. Third edition, 
\rm Ergebnisse der Math. und ihrer Grenzgebiete (2), Vol. {\bf 34}, Springer-Verlag, Berlin, 1994.

\Ref[MW] B. Moonen and T. Wedhorn, 
\sl Discrete invariants of varieties in positive characteristic,
\rm Int. Math. Res. Not. 2004, no. {\bf 72}, pp. 3855--3903.

\Ref[NV] 
M.-H. Nicole and A. Vasiu, 
\sl Traverso's isogeny conjecture for $p$-divisible groups,
\rm to appear in Rend. Semin. Mat. U. Padova, 8 pages.

\Ref[Oo1] 
F. Oort, 
\sl A stratification of a moduli space of abelian varieties, 
\rm Moduli of abelian varieties (Texel Island, 1999), pp. 345--416, Progr. of Math., Vol. {\bf 195}, Birkh\"auser, Basel, 2001.

\Ref[Oo2] 
F. Oort, 
\sl Foliations in moduli spaces of abelian varieties, 
\rm J. of Amer. Math. Soc. {\bf 17} (2004), no. 2, pp. 267--296.

\Ref[Oo3]
F. Oort,
\sl  Minimal $p$-divisible groups,  
\rm Ann. of Math. (2)  {\bf 161}  (2005),  no. 2, pp. 1021--1036. 

\Ref[Tr1]
C. Traverso,
\sl Sulla classificazione dei gruppi analitici di caratteristica positiva,
\rm Ann. Scuola Norm. Sup. Pisa {\bf 23}
 (1969), no. 3, pp. 481--507.

\Ref[Tr2]
C. Traverso,
\sl p-divisible groups over fields,
\rm  Symposia Mathematica, Vol. {\bf XI} (Convegno di Algebra Commutativa, INDAM, Rome, 1971), pp. 45--65, Academic Press, London, 1973.

\Ref[Tr3]
C. Traverso,
\sl Specializations of Barsotti--Tate groups,
\rm  Symposia Mathematica, Vol. {\bf XXIV} (Sympos., INDAM, Rome, 1979), pp. 1--21, Acad. Press, London-New York, 1981.  

\Ref[Va1]
A. Vasiu, 
\sl Integral canonical models of Shimura varieties of preabelian
 type, 
\rm Asian J. Math. {\bf 3} (1999), no. 2, pp. 401--518.

\Ref[Va2]
A. Vasiu,
\sl Crystalline boundedness principle,
\rm Ann. Sci. l'\'Ecole Norm. Sup. {\bf 39} (2006), no. 2, pp. 245--300.

\Ref[Va3]
A. Vasiu,
\sl Mod p classification of Shimura F-crystals,

\rm arxiv.org/abs/math/0304030.

\Ref[Va4]
A. Vasiu,
\sl Reconstructing p-divisible groups from their truncations of small level,
\rm arxiv.org/abs/math/0607268.

\Ref[Va5] 
A. Vasiu,
\sl Manin problems for Shimura varieties of Hodge type,

\rm arxiv.org/abs/math/0209410.

\Ref[Va6]
A. Vasiu,
\sl Integral canonical models of unitary Shimura varieties,

\rm arxiv.org/abs/math/0505507.

\Ref[We] T. Wedhorn, 
\sl The dimension of Oort strata of Shimura varieties of PEL-type,
\rm Moduli of abelian varieties (Texel Island, 1999), pp. 441--471, Progr. of Math., Vol. {\bf 195}, Birkh\"auser, Basel, 2001.

\Ref[Zi]
T. Zink,
\sl Isogenieklassen von Punkten von Shimuramannigfaltigkeiten mit Werten in einem endlichen K\"orper,
\rm Math. Nachr. {\bf 112} (1983), pp. 103--124.

}}
\noindent
\bigskip
\hbox{Adrian Vasiu}
\hbox{Department of Mathematical Sciences}
\hbox{Binghamton University}
\hbox{Binghamton, New York 13902-6000, U.S.A.}
\hbox{e-mail: adrian\@math.binghamton.edu}
\hbox{fax: 1-607-777-2450} 

\enddocument